\newif\ifalgin\algintrue 

\documentclass[12pt,twoside]{article}
\usepackage{latexsym}
\usepackage{graphicx}
\def\,{\mskip 3mu} \def\>{\mskip 4mu plus 2mu minus 4mu} \def\;{\mskip 5mu plus 5mu} \def\!{\mskip-3mu}
\def\dispmuskip{\thinmuskip= 3mu plus 0mu minus 2mu \medmuskip=  4mu plus 2mu minus 2mu \thickmuskip=5mu plus 5mu minus 2mu}
\def\textmuskip{\thinmuskip= 0mu                    \medmuskip=  1mu plus 1mu minus 1mu \thickmuskip=2mu plus 3mu minus 1mu}
\textmuskip
\def\beq{\dispmuskip\begin{equation}}    \def\eeq{\end{equation}\textmuskip}
\def\beqn{\dispmuskip\begin{displaymath}}\def\eeqn{\end{displaymath}\textmuskip}
\def\bqa{\dispmuskip\begin{eqnarray}}    \def\eqa{\end{eqnarray}\textmuskip}
\def\bqan{\dispmuskip\begin{eqnarray*}}  \def\eqan{\end{eqnarray*}\textmuskip}
\newenvironment{keywords}{\centerline{\bf\small
Keywords}\begin{quote}\small}{\par\end{quote}\vskip 1ex}
\def\paradot#1{\vspace{1.3ex plus 0.5ex minus 0.5ex}\noindent{\bf{#1.}}}
\def\paranodot#1{\vspace{1.3ex plus 0.5ex minus 0.5ex}\noindent{\bf{#1}}}
\newtheorem{theorem}{Theorem}

\def\aidx#1{}
\def\toinfty#1{\stackrel{#1\to\infty}{\longrightarrow}}
\def\eps{\varepsilon}
\def\epstr{\epsilon}                   
\def\nq{\hspace{-1em}}
\def\qed{\hspace*{\fill}$\Box\quad$\\}
\def\fr#1#2{{\textstyle{#1\over#2}}}
\def\SetR{I\!\!R}
\def\SetN{I\!\!N}
\def\SetB{I\!\!B}
\def\SetZ{Z\!\!\!Z}
\def\qmbox#1{{\quad\mbox{#1}\quad}}
\def\e{{\rm e}}                        
\def\B{\{0,1\}}
\def\v{\vec}
\def\es{\mbox{\o}}                     
\def\l{\ell}
\def\lb{{\log_2}}
\def\a{\alpha}
\def\b{\beta}
\def\G{\Gamma}                         
\def\Ga{\Gamma}                        
\def\Beta{\mbox{\rm Beta}}
\def\text#1{\mbox{\scriptsize #1}}
\def\EE{I\!\!E}

\begin{document}

\title{\vspace{-4ex}
\vskip 2mm\bf\Large\hrule height5pt \vskip 4mm
Exact Non-Parametric Bayesian Inference \\ on Infinite Trees
\vskip 4mm \hrule height2pt}
\author{{\bf Marcus Hutter}\\[3mm]
\normalsize RSISE$\,$@$\,$ANU and SML$\,$@$\,$NICTA \\
\normalsize Canberra, ACT, 0200, Australia%
\footnote{Preliminary results have been presented at the AISTATS 2005 conference \cite{Hutter:05bayestree}.}
\\
\normalsize \texttt{marcus@hutter1.net \ \  www.hutter1.net}
}
\date{31 March 2009}
\maketitle

\vspace{-4ex}\begin{abstract}
Given i.i.d.\ data from an unknown distribution, we consider the
problem of predicting future items. An adaptive way to estimate the
probability density is to recursively subdivide the domain to an
appropriate data-dependent granularity. In Bayesian inference one
assigns a data-independent prior probability to ``subdivide'', which
leads to a prior over infinite(ly many) trees. We derive an exact,
fast, and simple inference algorithm for such a prior, for the data
evidence, the predictive distribution, the effective model
dimension, moments, and other quantities. We prove asymptotic
convergence and consistency results, and illustrate the behavior of
our model on some prototypical functions.
\def\contentsname{\centering\normalsize Contents}
{\parskip=-2.7ex\tableofcontents}
\end{abstract}


\newpage
\section{Introduction}\label{secInt}

\paradot{Inference}
We consider the problem of inference from i.i.d.\ data $D$, in
particular of the unknown distribution $q$ the data is sampled
from. In case of a continuous domain this means inferring a
probability density from data. Without structural assumption on
$q$, this is hard to impossible, since a finite amount of data is
never sufficient to uniquely select a density (model) from an
infinite-dimensional space of densities (model class).

\begin{figure}
\centerline{\includegraphics[width=0.5\textwidth]{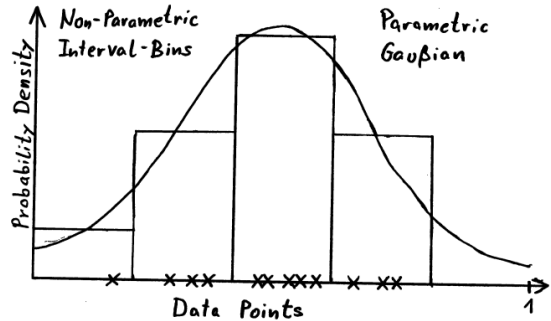}}
\caption{\label{figBinsGauss}Bins versus Gaussian estimate of
the true=data-generating probability density. (More decent diagrams
will be made for the final version).}
\end{figure}

\paradot{Methods}
In parametric estimation one assumes that $q$ belongs to a
finite-dimensional family. The two-dimensional family of Gaussians
characterized by mean and variance is prototypical (Figure
\ref{figBinsGauss}). The maximum likelihood (ML) estimate of $q$
is the distribution that maximizes the data likelihood. Maximum
likelihood overfits if the family is too large and especially if
it is infinite-dimensional. A remedy is to penalize complex
distributions by assigning a prior (2nd order) probability to the
densities $q$. Maximizing the model posterior (MAP), which is
proportional to likelihood times the prior, prevents overfitting.
A full Bayesian procedure keeps the complete posterior for
inference. Typically, summaries like the mean and variance of the
posterior are reported.

\paranodot{How to choose the prior?}
In finite or small compact low-dimensional spaces a uniform prior
often works (MAP reduces to ML). In the non-parametric case one
typically devises a hierarchy of finite-dimensional model classes of
increasing dimension. Selecting the dimension with maximal posterior
often works well due to the Bayes factor phenomenon
\cite{Good:83,Good:84,Jeffreys:35,Jaynes:03,MacKay:03}: In case the
true model is low-dimensional, higher-dimensional (complex) model
classes are automatically penalized, since they contain fewer
``good'' models. In a full Bayesian treatment one would assign a
prior probability (e.g.\ ${1\over d^2}$) to the dimension $d$  and
mix over the dimension.

\paradot{Interval Bins}
The probably simplest and oldest model for an interval domain is to
divide the interval (uniformly) into bins, assume a constant
distribution within each bin, and take a frequency estimate for the
probability in each bin (Figure \ref{figBinsGauss}), or a Dirichlet
posterior in Bayesian inference. There are heuristics for choosing
the number of bins as a function of the data size. The simplicity
and easy computability of the bin model is very appealing to
practitioners. Drawbacks are that distributions are discontinuous,
its restriction to one dimension (or at most low dimension: curse of
dimensionality), the uniform (or more generally fixed)
discretization, and the heuristic choice of the number of bins. We
present a full Bayesian solution to these problems, except for the
non-continuity problem. Our model can be regarded as an extension of
Polya trees \cite{Ferguson:73,Lavine:92,Lavine:94}.

\paradot{Related work}
There are plenty of alternative Bayesian models that
overcome some or all of the limitations. Examples are %
continuous Dirichlet process (mixtures) \cite{Ferguson:73}, %
Bernstein polynomials \cite{Petrone:02}, %
Bayesian field theory \cite{Lemm:03}, %
randomized Polya trees \cite{Paddock:03}, %
Bayesian bins with boundary averaging \cite{Endres:05}, %
Bayesian kernel density estimation %
or other mixture models \cite{Escobar:95},
and universal priors \cite{Hutter:04uaibook}, %
but exact analytical solutions are infeasible. %
Markov Chain Monte Carlo sampling \cite{Bishop:06}, %
Expectation Maximization algorithms \cite{Dempster:77}, %
variational methods \cite{Bishop:06}, %
efficient MAP or M(D)L approximations \cite{Kontkanen:07}, %
or kernel density estimation \cite{Gray:03} %
can often be used to obtain approximate numerical solutions, but
computation time and/or global convergence remain critical issues.
There are of course also plenty of non-Bayesian density estimators;
see (references in) \cite{Koller:98,Borovkov:98,Liu:07} in general,
and \cite{Kozlov:97,Koller:98} for density tree estimation in
particular.

\paradot{Our tree mixture model}
The idea of the model class discussed in this paper is very simple:
With some (e.g.\ equal) probability, we chose $q$ either uniform or
split the domain in two parts (of equal volume), and assign a prior
to each part, recursively, i.e.\ in each part again either uniform
or split. For finitely many splits, $q$ is a piecewise constant
function, for infinitely many splits it is virtually {\em any}
distribution. While the prior over $q$ is neutral about uniform
versus split, we will see that the posterior favors a split if and
only if the data clearly indicates non-uniformity. The method is a
full Bayesian non-heuristic tree approach to adaptive binning for
which we present a very simple and fast algorithm for computing
all(?) quantities of interest.

Note that we are not arguing that our model performs better in
practice than the more advanced models above. The main
distinguishing feature of our model is that it allows for a fast and
exact analytical solution. It's likely use is as a building block in
complex problems, where computation time and Bayesian integration
are the major issues. In any case, if/since the Polya tree model
deserves attention, also our model should.

\paradot{Contents}
In Section \ref{secTMM} we introduce our model and compare it to
Polya trees. We also discuss some example domains, like intervals,
strings, volumes, and classification tasks. %
Section \ref{secEPR} derives recursions for the posterior and the
data evidence. %
Section \ref{secAC} proves convergence/consistency. %
In Section \ref{secMQI} we introduce further quantities of interest,
including the effective model dimension, the tree size and height,
the cell volume, and moments, and present recursions for them. %
The proper case of infinite trees is discussed in Section
\ref{secIT}, where we analytically solve the infinite recursion at
the data separation level. %
Section \ref{secAlg} collects everything together and presents the
algorithm. %
In Section \ref{secEx} we numerically illustrate the behavior of
our model on some prototypical functions. %
Section \ref{secDisc} contains a brief summary, conclusions, and
outlook, including natural generalizations of our model. %
See \cite{Hutter:07btcode} for program code.

\section{The Tree Mixture Model}\label{secTMM}

\paradot{Setup and basic quantities of interest}
We are given i.i.d.\ data $D=(x^1,...,x^n)\in\G^n$ of size
$n$ from domain $\G$, e.g.\ $\G\subseteq\SetR^d\!$, sampled
from some unknown probability density $q:\G\to\SetR$. Standard
inference problems are to estimate $q$ from $D$ or to predict the
next data item $x^{n+1}\in\G$. By definition, the
(objective or aleatoric) data likelihood density under model $q$ is
\beq\label{eqLikelihood}
  \mbox{likelihood:}\qquad  p(D|q) \;\equiv\; q(x_1)\cdot...\cdot q(x_n)
\eeq
Note that we consider sorted data, which avoids annoying
multinomial coefficients. Otherwise this has no consequences.
Results are independent of the order and depend on the counts
only, as they should. %
A Bayesian assumes a (belief or $2^{nd}$-order or epistemic or
subjective) prior over models $q$ in some model class $Q$:
\beqn\label{eqPrior}
  \mbox{prior:}\qquad p(q) \qmbox{with} q\in Q
\eeqn
The data evidence is
\beq\label{eqEvidence}
  \mbox{evidence:}\qquad p(D) \;=\; \int_Q p(D|q)p(q) dq
\eeq
Having the evidence, Bayes' famous rule allows to compute
the (belief or $2^{nd}$-order or epistemic or subjective) posterior
of $q$:
\beq\label{eqPosterior}
  \mbox{posterior:}\qquad p(q|D) \;=\; {p(D|q)p(q)\over p(D)}
\eeq
The predictive distribution, i.e.\ the conditional probability that
next data item is $x=x^{n+1}$, given $D$, follows from the evidences
of $D$ and $(D,x)$:
\beq\label{eqPredDistr}
  \mbox{predictive distribution:}\qquad
  p(x|D) \;=\; {p(D,x)\over p(D)}
\eeq
Since the posterior is a complex object, we need summaries
like the expected $q$-probability of $x$ and (co)variances.
Fortunately they can also be reduced to computation of evidences:
\bqan\label{eqEqD}
  E[q(x)|D]
   &:=& \int q(x)p(q|D) dq
  \;=\; \int q(x){p(D|q)p(q)\over p(D)} dq
\\
   &=& {\int p(D,x|q)p(q)dq\over p(D)}
  \;=\; {p(D,x)\over p(D)}
  \;=\; p(x|D)
\eqan
where we used the formulas for the posterior, the
likelihood, the evidence, and the predictive distribution,
in this order. Similarly for the covariance we obtain
\bqan\label{eqCov}
  & & \mbox{Cov}[q(x)q(y)|D] \\
  & & \equiv\; E[q(x)q(y)|D]-E[q(x)|D]\!\cdot\!E[q(x)|D] \\
  & & =\; p(x,y|D)-p(x|D)p(y|D)
\eqan
We derive and discuss further summaries of $q$ for our particular tree
model, like the model complexity or effective dimension, the
tree height or cell size, and moments later.

\paradot{Hierarchical tree partitioning}
So far everything has been fairly general. We now introduce the
tree representation of domain $\G$. We partition $\G$ into
$\G_0$ and $\G_1$, i.e.\ $\G=\G_0\cup\G_1$ and
$\G_0\cap\G_1=\es$. Recursively we (sub)partition
$\G_z=\G_{z0}\dot\cup\G_{z1}$ for $z\in\SetB_0^m$,
where $\SetB_k^m:=\bigcup_{i=k}^m\{0,1\}^i$ is the set of all
binary strings of length between $k$ and $m$, and
$\G_\epstr=\G$, where $\epstr=\{0,1\}^0$ is the empty
string. We are interested in an infinite recursion, but for
convenience we assume a finite tree height $m<\infty$ and
consider $m\to\infty$ later. Also let $l:=\l(z)$ be the length of
string $z=z_1...z_l=:z_{1:l}$, and $|\G_z|$ the
volume or length or cardinality of $\G_z$.

\begin{figure}
\centerline{\includegraphics[width=0.9\textwidth]{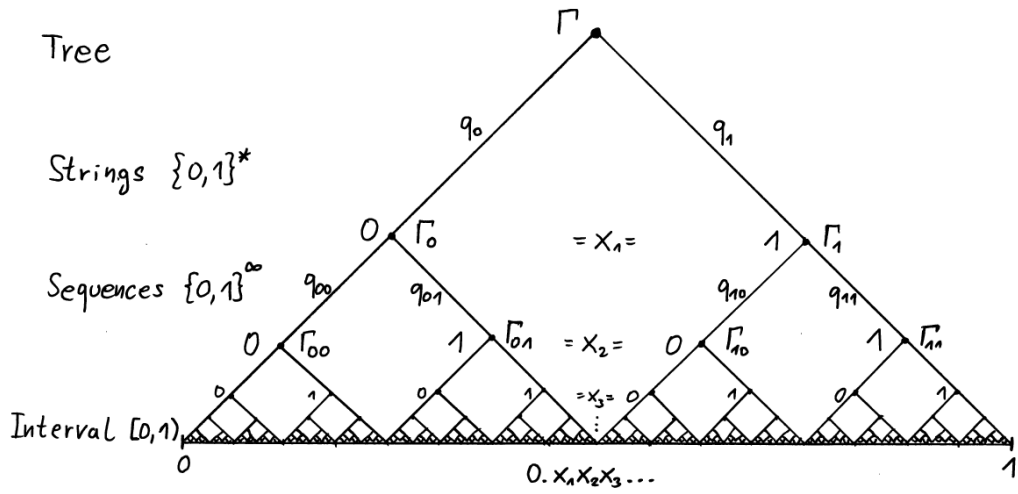}}
\caption{\label{figFullTree}Tree representation of strings or sequences or intervals.}
\end{figure}

\begin{figure}
\centerline{\includegraphics[width=0.3\textwidth]{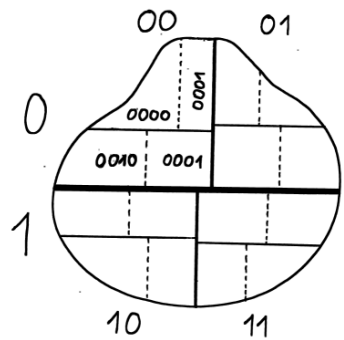}}
\caption{\label{figVolumeTree}Tree representation of volumes.}
\end{figure}

\paradot{Example spaces (Figures \ref{figFullTree} \& \ref{figVolumeTree})}
{\it Intervals:} Assume $\G=[0,1)$ is the unit interval, recursively
bisected into intervals $\G_z=[0.z, 0.z+2^{-l})$ of length
$|\G_z|=2^{-l}$, where $0.z$ is the real number in $[0,1)$
with binary expansion $z_1...z_l$.

{\it Strings:} Assume $\G_z=\{zy:y\in\B^{m-l}\}$ is the set of
strings of length $m$ starting with $z$. Then $\G=\B^m$ and
$|\G_z|=2^{m-l}$. For $m=\infty$ this set is continuous, for
$m<\infty$ finite.

{\it Trees:} Let $\G$ be a complete binary tree of height $m$ and
$\G_{z0}$ ($\G_{z1}$) be the left (right) subtree of
$\G_z$. If $|\G_z|$ is defined as one more than the number
of nodes in $\G_z$, then $|\G_z|=2^{m+1-l}$.

{\it Volumes:} Consider $\G\subset\SetR^d$, e.g.\ the
hypercube $\G=[0,1)^d$. We recursively halve $\G_z$ with a
hyperplane orthogonal to dimension $(l$ mod $d)+1$, i.e.\ we sweep
through all orthogonal directions. $|\G_z|=2^{-l}|\G|$.

{\it Compactification:} We can compactify
$\G\subseteq(1,\infty]$ (this includes
$\G=\SetN\setminus\{1\}$) to the unit interval
$\G':=\{{1\over x}:x\in\G\}\subseteq[0,1)$, and similarly
$\G\subseteq\SetR$ (this includes $\G=\SetZ$) to
$\G':=\{x\in[0,1):{2x-1\over x(1-x)}\in\G\}$.
All reasonable spaces
can be reduced to one of the spaces described above,
although this reduction may introduce unwanted artifacts.

{\it Classification:} Consider an observation $o\in\G'$ (e.g.\
email) that is classified as $c\in\{0,1\}$ (e.g.\ good versus
spam), where $\G'$ could be one of the spaces above (e.g.\
$o$ is a sequence of binary features in decreasing order of
importance). Then $x:=(o,c)\in\G:=\G'\times\{0,1\}$ and
$\G_{0z}=\G'_z\times\{0\}$ and
$\G_{1z}=\G'_z\times\{1\}$. Given $D$ (e.g.\
pre-classified emails), a new observation $o$ is classified as $c$
with probability $p(c|D,o)\propto p(D,x)$. Similar for more than
two classes.

In all these examples we have (chosen)
$|\G_{z0}|=|\G_{z1}|=\fr12|\G_z|$ $\forall
z\in\SetB_0^{m-1}$, and this is the only property we need and
henceforth assume. W.l.g.\ we assume/\linebreak[1]define/\linebreak[1]rescale $|\G|=1$.
Generalizations to non-binary and non-symmetric partitions are
straightforward and briefly discussed at the end.

\paradot{Identification}
We assume that $\{\G_z:z\in\SetB_0^m\}$ are (basis) events
that generate our $\sigma$-algebra. For every $x\in\G$ let
$x'$ be the string of length $\l(x')=m$ such that $x\in\G_{x'}$.
We assume that distributions $q$ are $\sigma$-measurable, i.e.\ to
be constant on $\G_{x'}$ $\forall x'\in\SetB^m$. For
$m=\infty$ this assumption is vacuous; we get {\em all} Borel
measures. Hence, we can identify the continuous sample space
$\G$ with the (for $m<\infty$ discrete) space $\SetB^m$ of
binary sequences of length $m$, i.e.\ in a sense all example
spaces are isomorphic.
While we have the volume model in mind for real-world
applications, the string model will be convenient for mathematical
notation, the tree metaphor will be convenient in discussion, and
the interval model will be easiest to implement and to present
graphically.

\paradot{Notation}
As described above, $\G$ may also be a tree. This
interpretation suggests the following scheme for defining the
probability of $q$ on the leaves $x'$. The probability of the left
child node $z0$, given we are in the parent node $z$, is
$P[\G_{z0}|\G_z,q]$, so we have
\beqn
   p(x|\G_z,q) = p(x|\G_{z0},q)\!\cdot\! P[\G_{z0}|\G_z,q]
   \qmbox{if} x\in\G_{z0}
\eeqn
and similarly for the right child. In the following we often have
to consider distributions conditioned to and in the subtree
$\G_z$, so the following notation will turn out
convenient\vspace{-1ex}
\beq\label{eqNotation}
  q_{z0} := P[\G_{z0}|\G_z,q], \quad
  q_{z1} := P[\G_{z1}|\G_z,q], \quad
  p_z(x|...) := 2^{-l}p(x|\G_z...)\vspace{-1ex}
\eeq
\beqn
  \Rightarrow p_z(x|q) = 2q_{zx_{l+1}} p_{zx_{l+1}}(x|q) =...
  = \!\!\prod_{i=l+1}^m\! 2q_{x_{1:i}}
  \;\mbox{if}\; x\in\G_z
\eeqn
where $p(x|\G_{x'},q):=|\G_{x'}|^{-1}=2^m$ is uniform (by
assumption). Note that $q_{z0}+q_{z1}=1$. Finally, let
\beqn
  \v q_{z*} \;:=\; (q_{zy}:y\in\SetB_1^{m-l})
\eeqn
be the ($2^{m-l+1}-2$)-dimensional {\em vector} or {\em ordered
set} or {\em tree} of all {\em reals} $q_{zy}\in[0,1]$ in subtree
$\G_z$. Note that $q_z\not\in\v q_{z*}$. The {\em
(non)density} $q_z(x):= p_z(x|q)$ depends on all and only these
$q_{zy}$. For $z\neq\epstr$, $q_z()$ and $p_z()$ are only
proportional to a density due to the factor $2^{-l}$, which has
been introduced to make $p_{x'}(x|...)\equiv 1$. (They are
densities w.r.t.\ $2^l\lambda_{|\G_z}$, where $\lambda$ is the
Lebesgue measure.) We have to keep this in mind in our
derivations, but can ignore this widely in our discussion.

\paradot{Polya trees}
In the Polya tree model one assumes that the $q_{z0}\equiv 1-q_{z1}$
are independent and Beta($\cdot,\cdot$) distributed, which defines
the prior over $q$. Polya trees form a conjugate prior class, since
the posterior is also a Polya tree, with empirical counts added to
the Beta parameters. If the same Beta is chosen in each node, the
posterior of $x$ is pathological for $m\to\infty$: %
The density does nowhere exist with probability 1. A cure is to
increase the Beta parameters with $l$, e.g.\ quadratically, but this
results in ``underfitting'' for large sample sizes, since
Beta(large,large) is too informative and strongly favors $q_{z0}$
near $\fr12$. It also violates scale invariance, which should
ideally hold if we do not have any prior knowledge about the scale.
That is, the p(oste)rior in $\G_0=[0,\fr12)$ should be the same as
for $\G=[0,1)$ (after rescaling all $x\leadsto x/2$ in $D$).

\paradot{The new tree mixture model}
The prior $p(q)$ follows from specifying a prior over $\v q_*$,
since $q(x)\propto q_{x_1}\cdot...\cdot q_{x_{1:m}}$ by
(\ref{eqNotation}). The distribution in each subset
$\G_z\subseteq\G$ shall be either $u$niform or non-uniform. A
necessary (but not sufficient) condition for uniformity is
$q_{z0}=q_{z1}=\fr12$.
\beq\label{tmmu1}
  p^u(q_{z0},q_{z1}) \;:=\; \delta(q_{z0}\!-\!\fr12)\delta(q_{z1}\!-\!\fr12),
\eeq
where $\delta()$ is the Dirac delta. To get uniformity on
$\G_z$ we have to recurse the tree down in this way.
\beq\label{tmmur}
  p_z^u(\v q_{z*}) \;:=\; p^u(q_{z0},q_{z1})p_{z0}^u(\v q_{z0*})p_{z1}^u(\v q_{z1*})
\eeq
with the natural recursion termination $p_z^u(\v q_{z*})=1$ when
$\l(z)=m$, since then $\v q_{z*}=\es$. For a non-uniform
distribution on $\G_z$ we allow any probability split
$q(\G_z)=q(\G_{z0})+q(\G_{z1})$, or equivalently
$1=q_{z0}+q_{z1}$. We assume a Beta prior on the $s$plit. Scale
invariance requires the Beta parameters to be the same in all
nodes of the tree and symmetry demands a symmetric Beta, i.e.\
\bqa\label{tmms1}
   p^s(q_{z0},q_{z1}) &:=&
   \Beta(q_{z0},q_{z1}|\a,\a),
\\ \label{defBeta}
  \nq\Beta(p,q|\a,\b) &:=&
  {\textstyle{\Ga(\a+\b)\over\Ga(\a)\Ga(\b)}}p^{\a-1}q^{\b-1}\delta(p\!+\!q\!-\!1)
\eqa
where $\Ga(\a)=\int_0^\infty t^{\a-1}\e^{-t}dt$ is the Gamma function.
For $\a=1$ this specializes to the natural uniform prior
$p^s(q_{z0},q_{z1})=\delta(q_{z0}+q_{z1}-1)$ on the $s$plit.
We now recurse down the tree
\beq\label{tmmsr}
   p_z^s(\v q_{z*}) \;:=\; p^s(q_{z0},q_{z1})p_{z0}(\v q_{z0*})p_{z1}(\v q_{z1*})
\eeq
again with the natural recursion termination $p_z(\v q_{z*})=p(\es)=1$
when $\l(z)=m$. Finally we have to mix the uniform with the
non-uniform case.
\beq\label{tmmmix}
  p_z(\v q_{z*}) \;:=\; u\!\cdot\!p_z^u(\v q_{z*})+s\!\cdot\!p_z^s(\v q_{z*})
\eeq
with $u,s\in[0,1]$ and $u+s=1$. The 50/50 mixture $u=s=\fr12$ will
be of special interest. This completes the specification of the
prior $p(q)=p_\epstr(\v q_*)$.%
\footnote{Note that $p_z(\v q_{z*})$ is {\em not} the
marginal of $p(q)$ to $\v q_{z*}$, but one can show that
$p_z(\v q_{z*})=p(\v q_{z*}|q_{z_1}\neq\fr12,...,q_{z_{1:l}}\neq\fr12)$
and optionally additional conditions on some or all $q\not\in\v q_{z*}$.}

For example, if the first bit in $x$ is a class label
and the remaining are binary features in decreasing order of
importance, then given class and features $z=x_{1:l}$, further features
$x_{l+1:m}$ could be relevant for classification ($q_z(x)$ is
non-uniform) or irrelevant ($q_z(x)$ is uniform).

\paradot{Comparison to the Polya tree}
Note the important difference in the recursions (\ref{tmmur}) and
(\ref{tmmsr}). Once we decided on a uniform distribution
(\ref{tmmu1}) we have to equally split probabilities down the
recursion to the end, i.e.\ we recurse in (\ref{tmmur}) with $p^u$,
rather than the mixture $p$ (this actually allows to solve the
recursion). On the other hand if we decided on a non-uniform split
(\ref{tmms1}), the left and right partition each itself may be
uniform or not, i.e.\ we recurse in (\ref{tmmsr}) with the mixture
$p$, rather than $p^s$. Inserting (\ref{tmms1}) in (\ref{tmmsr}) in
(\ref{tmmmix}) and recursively (\ref{tmmu1}) in (\ref{tmmur}) in
(\ref{tmmmix}) we get the following recursion for the prior
\beq\label{EqPqRec}
  p_z(\v q_{z*}) = u\!\cdot\!\nq\prod_{y\in\SetB_1^{m-l}}\nq\; \delta(q_{zy}\!-\!\fr12)
  + s\!\cdot\!\Beta(q_{z0},q_{z1}|\a,\a)p_{z0}(\v q_{z0*})p_{z1}(\v q_{z1*})
\eeq
Choosing $u=0$ would lead to the Polya tree model (and its
problems) with $q_{z0}\sim$ Beta$(\cdot|\a,\a)$.
With $p$ instead of $p^u$ on the r.h.s.\ of
(\ref{tmmur}) we would get a quasi-Polya model (same problems) with
$q_{z0}\sim u\cdot\mbox{Beta}(\cdot|\infty,\infty)+s\cdot\mbox{Beta}(\cdot|\a,\a)$.

For $m\to\infty$, our model is ``scale'' invariant {\em and} leads
to continuous densities for $n\to\infty$, unlike the Polya tree
model. We also don't have to tune Beta parameters. We can use a
non-informative prior like $\a=1$ and $u=s=\fr12$. The model ``tunes
itself'' by suitably assigning high/low posterior probability to
subdividing cells. While Polya trees form a natural conjugate prior
class, our prior does not directly, but may be generalized to do so.
The computational complexity for the quantities of interest will be
the same (essentially $O(n)$), i.e.\ as good as it could be.

\paradot{Formal and effective dimension}
Formally our model is $2\cdot(2^m-1)$-dimensional,
but the effective dimension can by much smaller, since $\v q_*$ is
forced with a non-zero probability to a much smaller polytope, for
instance with probability $u$ to the zero-dimensional globally
uniform distribution. We will compute the effective p(oste)rior
dimension. Alternatively, we could have considered a mixture over
all ($\widehat=$ lower dimensional) partial trees with $\G_z$ as
leaf if $q$ is uniform on $\G_z$, but considering one complete tree
is more convenient for analytical manipulation.

\section{Evidence and Posterior Recursion}\label{secEPR}

At the end of Section \ref{secTMM} we defined our tree mixture
model. The next step is to compute the standard quantities of
interest defined at the beginning of Section \ref{secTMM}. The
evidence (\ref{eqEvidence}) is key, the other quantities (posterior,
predictive distribution, expected $q(x)$ and its variance) follow
then immediately. Let
\beqn
  D_z \;:=\; \{x\in D : x\in\G_z\}
\eeqn
be the $n_z:=|D_z|$ data points that lie in subtree $\G_z$.
We compute $p_z(D_z)$ recursively for all $z\in\SetB_0^{m-1}$, which gives
$p(D)=p_\epstr(D_\epstr)$.

\begin{theorem}[Evidence recursion]\label{thEvRec}
For $\l(z)<m$ the recursion for the evidence is
\bqa\label{eqEvDens}\label{eqEvRec}
  p_z(D_z) &=& u \;+\; s\!\cdot\!{p_{z0}(D_{z0})p_{z1}(D_{z1})\over w({n_{z0},n_{z1}})}
\\ \label{eqWeights}
     w({n_{z0},n_{z1}})
   \;\nq &:=& \nq\; {2^{-n_z}\!\cdot\! \Ga(n_z\!+\!2\a)\over \Ga(n_{z0}\!+\!\a)\G(n_{z1}\!+\!\a)}
       \!\cdot\! {\Ga(\a)^2\over\Ga(2\a)}
   =: w_{n_z}(\Delta_z)
\\ \nonumber
  \qquad n_z &=& n_{z0}+n_{z1}, \quad
  \Delta_z \;:=\; {n_{z0}\over n_z}-\fr12
\eqa
The recursion terminates with $p_z(D_z)\equiv 1$ when $\l(z)=m$.
\end{theorem}

\noindent The recursion (\ref{eqEvRec}) follows by multiplying
(\ref{eqg}) in Theorem \ref{thprec} (stated and proven below) with
$p_z(D_z)$ and adding $u$.
For $\a=1$, (\ref{eqEvDens}) and in particular the
weight $w_{n_z}=2^{-n_z}(n_z+1)({n_z\atop n_{z0}})$ can be
interpreted as follows: With probability $u$, the evidence is
uniform in $\G_z$. Otherwise data $D_z$ is split into two
partitions of size $n_{z0}$ and $n_{z1}=n_z-n_{z0}$. First, choose
$n_{z0}$ uniformly in $\{0,...,n_z\}$. Second, given $n_z$, choose
uniformly among the $({n_z\atop n_{z0}})$ possibilities of
selecting $n_{z0}$ out of $n_z$ data points for $\G_{z0}$ (the
remaining $n_{z1}$ are then in $\G_{z1}$). Third, distribute
$D_{z0}$ according to $p_{z0}(D_{z0})$ and $D_{z1}$ according to
$p_{z1}(D_{z1})$. Then, the evidence in case of a split is the
second term in (\ref{eqEvRec}). The factor $2^{n_z}$ is due to our
normalization convention (\ref{eqNotation}). This also verifies
that the r.h.s.\ yields the l.h.s.\ if integrated over all $D_z$,
as it should be.
%
For $n_z\to\infty$ we will show in Section \ref{secAC} that
$w_{n_z}\to\infty$ if $n_{z0}\sim n_{z1}\to\infty$ and
$w_{n_z}\to 0$ otherwise, indicating that the weight $w$ is large
(small) for (non)uniform distribution, as it should be.

\begin{theorem}[Posterior recursion]\label{thprec}
For $\l(z)<m$ the recursion for the posterior is
\bqa\label{eqprec}
  \hspace*{-10ex} p_z(\v q_{z*}|D_z) & = &
  {u\over p_z(D_z)} \prod_{y\in\SetB_1^{m-l}}\!\!\! \delta(q_{zy}-\fr12)
\\ \nonumber
  & + & g_z(D_z) \Beta(q_{z0},q_{z1}|n_{z0}\!+\!\a,n_{z1}\!+\!\a)
   p_{z0}(\v q_{z0*}|D_{z0}) p_{z1}(\v q_{z1*}|D_{z1})
\\ \label{defg}\label{eqg}\label{eqSplitProb}
   g_z(D_z) &:=& s\!\cdot\!{p_{z0}(D_{z0})p_{z1}(D_{z1})
                   \over p_z(D_z) w({n_{z0},n_{z1}}) }
  \;=\; 1-{u\over p_z(D_z)}
\eqa
The recursion terminates with $p_z(\v q_{z*}|D_z)\equiv 1$ when $\l(z)=m$.
\end{theorem}

$g_z(D_z)$ may be interpreted as the posterior probability of
splitting $\G_z$.

\paradot{Proof}
Using Bayes rule (\ref{eqPosterior}) we represent the posterior as
\beq\label{eqBqD}
  p_z(\v q_{z*}|D_z)p_z(D_z) \;=\; p_z(D_z|\v q_{z*})p_z(\v q_{z*})
\eeq
and further substitute $p_z(\v q_{z*})= u p_z^u(\v q_{z*})+s p_z^s(\v q_{z*})$ (\ref{tmmmix}).
For the uniform part we get
\bqa\nonumber
  p_z(D_z|\v q_{z*})\!\cdot\!p_z^u(\v q_{z*})
  &=& \nq\;\;\prod_{x\in D_z}(2q_{x_{1:l+1}}\!\cdot\!...\!\cdot\!2q_{x_{1:m}})
      \!\cdot\nq\prod_{y\in\SetB_1^{m-l}}\nq \delta(q_{zy}-\fr12)
\\ \label{eqERPa}\label{eqERPb}\label{eqERPb1}
  &=& \nq\prod_{y\in\SetB_1^{m-l}}\nq \delta(q_{zy}-\fr12),
\eqa
where we recursively inserted (\ref{tmmu1}) in (\ref{tmmur}), and
(\ref{eqLikelihood}) and (\ref{eqNotation}) into (\ref{eqERPb1}).
Due to the $\delta$, we can simply set all $q_{zy}=\fr12$.
For the split we get
\bqa\nonumber
   & & \nq\nq p_z(D_z|\v q_{z*})\!\cdot\!p_z^s(\v q_{z*})
\\ \label{eqERPb2}
   &\nq=& \Big(\!\prod_{x\in D_{z0}\nq} 2q_{z0}\!\Big) p_{z0}(D_{z0}|\v q_{z0*})
          \Big(\!\prod_{x\in D_{z1}\nq} 2q_{z1}\!\Big) p_{z1}(D_{z1}|\v q_{z1*})
\\ \nonumber
   & & \times\; {\textstyle{\Ga(2\a)\over\Ga(\a)^2}}\, q_{z0}^{\a-1}q_{z1}^{\a-1}
             \delta(q_{z0}\!+\!q_{z1}\!-\!1)p_{z0}(\v q_{z0*})p_{z1}(\v q_{z1*})
\\ \label{eqERPc}
  &\nq=& 2^{n_z}{\textstyle{\Ga(2\a)\over\Ga(\a)^2}}
            \,q_{z0}^{n_{z0}+\a-1}q_{z1}^{n_{z1}+\a-1}\delta(q_{z0}\!+\!q_{z1}\!-\!1)
\\ \nonumber & & \times p_{z0}(\v q_{z0*}|D_{z0})p(D_{z0})
                        p_{z1}(\v q_{z1*}|D_{z1})p(D_{z1})
\\ \label{eqERPd}
  &\nq=& {\textstyle{1\over s}}g_z(D_z)p_z(D_z)
  \Beta(q_{z0},q_{z1}|n_{z0}+\a,n_{z1}+\a)
\\ \nonumber
  & & \times p_{z0}(\v q_{z0*}|D_{z0}) p_{z1}(\v q_{z1*}|D_{z1})
\eqa
In (\ref{eqERPb2}) we split $D_z$ into $D_{z0}$ and $D_{z1}$ and
used (\ref{eqLikelihood}) and (\ref{eqNotation}) and
the fact that $p_{z0}(D_{z0}|\v q_{z*})$ depends on $q$
through $\v q_{z0*}$ only. We also inserted (\ref{defBeta}) in
(\ref{tmms1}) in (\ref{tmmsr}) in (\ref{eqERPb2}) and used
$n_{z0}+n_{z1}=n_z$. Rearranging terms and using Bayes rule
(\ref{eqBqD}) for subtrees $\G_{z0}$ and $\G_{z1}$ we get
(\ref{eqERPc}). The last equality is easiest proven backwards by
inserting $g_z$ (\ref{defg}) and $w$ (\ref{eqWeights}) and \Beta\
(\ref{defBeta}) into (\ref{eqERPd}). Inserting (\ref{tmmmix}) and
(\ref{eqERPa}) and (\ref{eqERPb2})-(\ref{eqERPd}) into
(\ref{eqBqD}) and dividing by $p_z(D_z)$ yields (\ref{eqprec}).

Integrating (\ref{eqprec}) over $\v q_{z*}$ and noting that $\int d\v
q_{z*}=\int dq_{z0}dq_{z1}\cdot\int d\v q_{z0*}\cdot\int d\v q_{z1*}$
factorizes and that $\prod\delta()$ and \Beta() and the $\v
q_{z*}$, $\v q_{z0*}$ and $\v q_{z1*}$ posteriors are all
proper densities which integrate to 1, we get
\beqn
  1 \;=\; {u\over p_z(D_z)}\!\cdot\!1
          + g_z(D_z)\!\cdot\!1\!\cdot\!1\!\cdot\!1
\eeqn
This shows the last equality in (\ref{eqg}). Theorem
\ref{thEvRec} (\ref{eqEvRec}) now follows by multiplying
(\ref{eqg}) with $p_z(D_z)$ and adding $u$.

For a formal proof of the recursion termination, recall
(\ref{eqNotation}): For $\l(z)=m$ and $x\in\G_z$ we have
$\G_{x'}=\G_z$ $\Rightarrow$ $p_z(x|q)=1$ $\Rightarrow$
$p_z(D_z|q)=1$ $\Rightarrow$ $p_z(D_z)=1$. \qed

\section{\boldmath Asymptotic Convergence/Consistency ($n\to\infty$)}\label{secAC}

\paradot{\boldmath Discussing the weight}
The relative probability of splitting (second term on r.h.s.\ of
(\ref{eqEvDens})) to the uniform case (first term in r.h.s.\ of
(\ref{eqEvDens})) is controlled by the weight $w$. Large (small)
weight indicates a (non)uniform distribution, provided $p_{z0}$
and $p_{z1}$ are $O(1)$. The balance $\Delta_z\approx 0$
($\not\approx 0$) indicates a (non)symmetric partitioning of the
data among the left and right branch of $\G_z$. Asymptotically for
large $n_z$ (and small $\Delta_z$), we have
\beqn
  w_{n_z}(\Delta_z) \approx c_\alpha{\textstyle\sqrt{2n_z\over\pi}}\;\e^{-2n_z\Delta_z^2}
\eeqn
where $c_\a>0$ is some finite constant.
Assume that data $D$ is sampled from the true distribution $\dot
q$. The probability of the left branch $\G_{z0}$ of $\G_z$
is $\dot q_{z0}\equiv P[\G_{z0}|\G_z,\dot q]=
2^l\dot q_z(\G_{z0})$. The relative frequencies ${n_{z0}\over
n_z}$ asymptotically converge to $\dot q_{z0}$. More precisely
${n_{z0}\over n_z}=\dot q_{z0}\pm O(n_z^{\smash{-1/2}})$.
Similarly for the right branch.
Assume the probabilities are equal ($\dot
q_{z0}=\dot q_{z1}=\fr12$), possibly but not necessarily due to a
uniform $\dot q_z()$ on $\G_z$. Then $\Delta_z=O(n_z^{\smash{-1/2}})$,
which implies
\beqn\label{eqwtoinfty}
  w_{n_z}(\Delta_z) \approx \Theta(\sqrt{n}_z)
  \;\toinfty{n_z}\; \infty
  \qmbox{if} \dot q_{z0}=\dot q_{z1}=\fr12,
\eeqn
consistent with our anticipation.
Conversely, for $\dot q_{z0}\neq\dot q_{z1}$ (which implies
non-uniformity of $\dot q_z()$) we have $\Delta_z\to c:=\dot
q_{z0}-\fr12 \neq 0$, which implies
\beqn\label{eqwto0}
  w_{n_z}(\Delta_z) \approx {\textstyle\sqrt{2n_z\over\pi}}\,\e^{-2n_z c^2}
  \;\toinfty{n_z}\; 0
  \qmbox{if} \dot q_{z0}\neq\dot q_{z1},
\eeqn
again, consistent with our anticipation.
Formally, the following can be proven:

\begin{theorem}[Weight asymptotics]\label{wAsym}
For $\dot q_{z0}=\fr12$ we have with probability 1 (w.p.1)
\beqn
  i)   \quad \lim_{n_z\to\infty}{\ln n_z\over\sqrt{n_z}}\,w_{n_z}(\Delta_z) = \infty, \quad\mbox{and} \\
\eeqn
\beqn
  ii)  \quad \mathop{\lim\sup}_{n_z\to\infty}\sqrt{\pi\over 2n_z\!\!}\;w_{n_z}(\Delta_z)
       = c_\alpha \;\;\{ {\textstyle{>\,0\atop <\infty}}
\eeqn
where $c_\a=4^{\a-1}\Ga(\a)^2/\Ga(2\a)$.
For $\dot q_{z0}\neq\fr12$ we have w.p.1
\beqn
  iii) \quad \lim_{n_z\to\infty}\e^{2n_zc^2}w_{n_z}(\Delta_z)=0 \quad\forall\; |c|<|\dot q_{z0}\!-\!\fr12|
\eeqn
\end{theorem}

\paradot{Proof}
We will drop the index $z$ everywhere. We need an asymptotic
representation of $w$ for $n_0,n_1\to\infty$. Using Stirling's
approximation $\ln\Ga(x)=(x-\fr12)\ln x-x+\fr12\ln(2\pi)+O({1\over
x})$ we get after some algebra
\bqa\label{lnw}
  \ln w_n(\Delta) &=&
      -n[H(\fr12)\!-\!H(\fr12\pm\tilde\Delta)]
      + \fr12\ln{\textstyle{n\over 2\pi}}
\\ \nonumber
      & & \nq+\,(2\a\!-\!1)H(\fr12\pm\tilde\Delta) +C_\a + \textstyle O({\a\over n_0}\!+\!{\a\over n_1}),
\\ \nonumber
  H(p) &=& -p\ln p -(1\!-\!p)\ln(1\!-\!p) \;=\; \mbox{Entropy}(p),
\\ \nonumber
  \tilde\Delta &=& \textstyle {n\over n+2\a-1}\Delta_n
  \;=\; {{1\over 2}(n_0-n_1)\over n+2\a-1},
  \quad \Delta_n=\Delta={\textstyle{n_0\over n}}-\fr12,
\\ \nonumber
  C_\a &=& 2\ln\Ga(\a)-\ln\Ga(2\a)
\eqa

{\bf (i)} follows from the law of the iterated logarithm
\beqn
  \mathop{\lim\sup}_{n\to\infty}
  {|X_1+...+X_n-n\mu|\over\sigma\sqrt{n\ln\ln n}} \;=\;1 \quad\mbox{w.p.1}
\eeqn
for i.i.d.\ random variables $X_1,...,X_n$ with mean $\mu$ and
variance $\sigma^2$. For the $i^{th}$ data item in $D$, let
$X_i=1$ if $x\in D_0$ and $X_i=0$ if $x\in D_1$. Then the $X_i$
are i.i.d.\ Bernoulli($\dot q_0$) with $\mu=\dot q_0=\fr12$ and
$\sigma^2=\dot q_0\dot q_1={1\over 4}$. Further, $X_1+...+X_n=n_0$ implies
$X_1+...+X_n-n\mu=n\Delta_n$ implies
$\lim\,\sup_n\sqrt{4n\over\ln\ln n}|\Delta_n|=1$ w.p.1.\ implies
\beqn
   \Delta_n^2 \;\leq\; (1\!+\!\eps){\textstyle{\ln\ln n\over 4n}} \quad\mbox{w.p.1}
\eeqn
for all sufficiently large $n$ and any $\eps>0$. Using
$\tilde\Delta=\Delta+O(\fr1n)$, (\ref{lnw}) can be further
approximated by
\beqn
  \ln w_n(\Delta) = -n[H(\fr12)-H(\fr12\!\pm\!\Delta)] + \fr12\ln{\textstyle{n\over 2\pi}} + O(1)
\eeqn
A Taylor series expansion around $\Delta=0$ yields
\beqn
  H(\fr12)-H(\fr12\!\pm\!\Delta) \;=\; 2\Delta^2+O(\Delta^4)
  \;\leq\; (1\!+\!\eps){\textstyle{\ln\ln n\over 2n}} + O(({\textstyle{\ln\ln n\over 4n}})^2)
\eeqn
which implies
\beqn
  \ln w_n(\Delta_n) - \fr12\ln{\textstyle{n\over 2\pi}} + \ln\ln n
  \;\geq\; \fr12(1-\eps)\ln\ln n +O(1) \;\longrightarrow\; \infty \qmbox{for} \eps<1
\eeqn
which implies $(i)$ by exponentiation.

{\bf(ii)} $(a)$ Convexity of $\ln\Ga(x)$ implies that $\ln
w_n(\Delta)$ is concave and symmetric in $\Delta$, hence
$\ln w_n(\Delta)$ assumes its global maximum at $\Delta=0$. $(b)$ From
(\ref{lnw}) it follows that $\ln w_n(0)=\fr12\ln{n\over
2\pi}+(2\a-1)H(\fr12) + C_\a + O(\fr1n)$. $(c)$
$(2n\Delta_n)_{n=1}^\infty$ is a symmetric random walk, hence
infinitely often passes zero w.p.1. $(a)$ and $(b)$ imply the $\leq$ and
$(b)$ and ($c)$ the $\geq$ in $\lim\,\sup_n[\ln
w_n(\Delta_n)-\fr12\ln{n\over 2\pi}] = (2\a-1)\ln 2 + C_\a$ w.p.1.
Exponentiation yields $(ii)$.

{\bf(iii)} Since $\Delta_n\to\dot q_0-\fr12$ w.p.1, (\ref{lnw})
implies $2nc^2+\ln w_n\sim n[2c^2-H(\fr12)+H(\dot q_0)]\to
-\infty$, since $H(\fr12)-H(\dot q_0)\geq 2(\dot q_0-\fr12)^2>2c^2$.
The asymptotic representation also holds for $n_0=0$ or $n_1=0$,
hence $(iii)$ follows by exponentiation for all $\dot q_0$,
including 0 and 1. \qed

\paradot{\boldmath Asymptotics of the evidence $p(D)$}
The typical use of the posterior $p(x|D)$ is as an estimate for
the unknown true distribution $\dot q(x)$. This makes sense if $p(x|D)$
is close to $\dot q(x)$. We show that the finite tree mixture model is
indeed consistent in the sense that $p(x|D)$
converges\footnote{All $\toinfty{n}$ statements hold with
probability 1 (w.p.1).} to $\dot q(x)$ and the posterior of $q()$
concentrates around the true value $\dot q()$ for $n\to\infty$.

\begin{theorem}[Evidence asymptotics]\label{thEvAs}
For fixed $m<\infty$ and $n_z\to\infty$, the posterior
$p_z(x|D_z)\to\dot q_z(x)$ for all $x\in\G_z$. Furthermore, for the
evidence w.p.1 we have
\beqn
  p_z(D_z)\;\left\{\begin{array}{ccl}
    \stackrel{poly.}\longrightarrow & u & \mbox{for uniform $\dot q_z()$ and $l<m$,} \\
    \equiv                         & 1 & \mbox{for $l=m$,} \\
    \stackrel{exp.}\longrightarrow & \infty & \mbox{for non-uniform $\dot q_z()$ provided $s>0$.} \\
  \end{array} \right.
\eeqn
\end{theorem}

\paranodot{Proof}
by induction on $l$. We have to show slightly more, namely also that
$p_z(D_z,x)\to c\in\{u,1,\infty\}$. For $l=m$, the theorem is
obvious, since $\dot q_z(x)$ must be uniform on $\G_z$ and
$p_z(D_z)\equiv 1\equiv p_z(D_z,x)$, hence $p_z(x|D_z)\equiv 1\equiv
\dot q_z(x)$.
Now assume the theorem holds for $\G_{z0}$ and $\G_{z1}$
and $l<m$. We show that it then also holds for $\G_z$.
Assume $u>0$ first.

$(a)$ Assume first, that $\dot q_z()$ is uniform. This implies
that also $\dot q_{z0}()$ and $\dot q_{z1}()$ are uniform, hence
$n_{z0},n_{z1}\to\infty$, hence by induction hypothesis,
$p_{z0}(D_{z0})$ and $p_{z1}(D_{z1})$ are bounded. Further,
$w_{n_z}(\Delta_z)\stackrel{poly.}\longrightarrow\infty$ for
$n_z\to\infty$ (by Theorem \ref{wAsym}$(i)$ and $(ii)$). Hence,
$p_z(D_z)\stackrel{poly.}\longrightarrow u$ from (\ref{eqEvDens}).
Similarly $p_z(D_z,x)\to u$, hence $p_z(x|D_z)\to 1\equiv\dot
q_z(x)$ for $x\in\Gamma_z$.

$(b)$ We now consider the case of non-uniform $\dot q_z()$.
$(i)$ Consider the case that $\dot q_{z0}()$ or $\dot q_{z1}()$
(or both) are non-uniform first. $p_{z0}(D_{z0})\geq u>0$ and
$p_{z1}(D_{z1})\geq u>0$, and one of them diverges exponentially.
Since $w_{n_z}$ grows at most with $O(\sqrt{n}_z)$ by Theorem
\ref{wAsym}$(ii)$, we see from (\ref{eqEvDens}) that also
$p_z(D_z)\sim s\cdot p_{z0}(D_{z0})p_{z1}(D_{z1})/w_{n_z}$ diverges
exponentially, and similarly $p_z(D_z,x)$.
$(ii)$ If both $\dot q_{z0}()$ and $\dot q_{z1}()$ are uniform,
then $q_{z0}\neq\fr12$, since we assumed non-uniform $\dot q_z()$.
This implies bounded $p_{z0}(D_{z0})$ and $p_{z1}(D_{z1})$,
but exponentially diverging $w_n^{-1}$ by Theorem
\ref{wAsym}$(iii)$. Hence, again, $p_z(D_z)$ and similarly
$p_z(D_z,x)$ diverge exponentially. In both cases, $(i)$ and
$(ii)$, assuming w.l.g.\ $x\in\G_{z0}$, the ratio is
\beq\label{eqEvAsPr1}\textstyle
  p_z(x|D) \sim
   {w_{n_z}\over w_{n_z+1}}\!\cdot\!p_{z0}(x|D_{z0})
     = 2\!\cdot\!{n_{z0}+\a\over n_z+2\a}\!\cdot\!p_{z0}(x|D_{z0})
  \sim 2\!\cdot\! \dot q_{z0}\!\cdot\!\dot q_{z0}(x) = \dot q_z(x)
\eeq
See (\ref{eqNotation}) for notation and how the density factor 2
disappears.
For $u=0$, (\ref{eqEvAsPr1}) holds for any $\dot q_z()$. Further,
$p_z(D_z)\stackrel{exp.}\longrightarrow\infty$ still holds, since
$p_{z0}(D_{z0})$ tends not faster than polynomially to zero by
induction.
\qed

\begin{theorem}[Posterior consistency]\label{thpc}
The posterior of $\v q_{z*}$ concentrates for $n_z\to\infty$ around
the true value $\,\dot{\!\v q}_{z*}$ w.p.1., i.e.\footnote{The topology of
weak convergence or convergence in distribution is used.}
\beq\label{eqpc}
  p_z(\v q_{z*}|D_z) \quad
  \mathop{\longrightarrow}\limits^{n_z\to\infty}_{w.p.1}
  \prod_{y\in\SetB_1^{m-l}} \delta(q_{zy}-\dot q_{zy})
\eeq
\end{theorem}

\paradot{Proof}
We prove consistency (\ref{eqpc}) by induction over $l$. For
$\l(z)=m$ the l.h.s.\ and r.h.s.\ are formally 1, since a density
over an empty space and an empty product are defined as 1. Assume
that consistency (\ref{eqpc}) holds for $\l(z0)=\l(z1)=l+1$. For
$n_z\to\infty$, the \Beta\ concentrates around ${n_{z0}\over
n_z}\to\dot q_{z0}$ and ${n_{z1}\over n_z}\to\dot q_{z1}$ w.p.1:
\beqn
   \Beta(q_{z0},q_{z1}|n_{z0}\!+\!\a,n_{z1}\!+\!\a)
   \to \delta(q_{z0}\!-\!\dot q_{z0})\delta(q_{z1}\!-\!\dot q_{z1})
\eeqn
Inserting this and (\ref{eqpc}) for $z0$ and $z1$ into the r.h.s.\
of recursion (\ref{eqprec}) we get
\beq\label{eqpcr}
  p_z(\v q_{z*}|D_z) \to
  {u\over p_z(D_z)} \prod_{y\in\SetB_1^{m-l}}\!\!\! \delta(q_{zy}-\fr12)
  + \Big(1-{u\over p_z(D_z)}\Big) \prod_{y\in\SetB_1^{m-l}}\!\!\! \delta(q_{zy}-\dot q_{zy})
\eeq
For uniform $\,\dot{\!\v q}_{z*}$, i.e.\ $\fr12=\dot q_{zy}$
$\forall y\in\SetB_1^{m-l}$ the r.h.s.\ reduces to the r.h.s.\ of
(\ref{eqpc}). For non-uniform $\,\dot{\!\v q}_{z*}$, Theorem
\ref{thEvAs}$(iii)$ shows that $p_z(D_z)\to\infty$ (exponentially), and
the r.h.s.\ of (\ref{eqpcr}) converges (rapidly) to the r.h.s. of (\ref{eqpc}).
\qed

\section{More Quantities of Interest}\label{secMQI}

In this section we introduce further quantities of interest and
present recursions for them. They all can be written as expectations
\beq\label{eqfexp}
  E_z[f(\v q_{z*})|D_z]
  \;:=\; \int f(\v q_{z*}) p_z(\v q_{z*}|D_z)d\v q_{z*}
\eeq
for suitable functions $f$ (and similarly for
$P_z[...|...]=P[...|\G_z...]$). For instance for the evidence
(\ref{eqEvDens}) we used $f(\v q_{z*})\equiv 1$ and for the
posterior (\ref{eqprec}) formally $f(\v
q_{z*})=\prod_y\delta(q_{zy}-q'_{zy})$. Below we consider the
model dimension, cell number, tree height, cell size, and moments. The
derivations of the recursion all follow the same scheme, inserting
the (recursive) definition of $f$ and the recursion
(\ref{eqprec}) into the r.h.s.\ of (\ref{eqfexp}), and
rearranging and identifying terms. These details will be omitted.

\paradot{Model dimension and cell number}
As discussed in Section \ref{secTMM}, the effective dimension of
$\v q_*$ is the number of components that are not forced to $\fr12$
by (\ref{tmmu1}). Note that a component may be ``accidentally''
$\fr12$ in (\ref{tmms1}), but since this is an event of probability
0, we don't have to care about this subtlety. So the effective
dimension $N_{\v q_{z*}}=\#\{q\in\v q_{z*}:q\neq\fr12\}$ of $\v
q_{z*}$ can be given recursively as
\beq\label{eqNrec}
   N_{\v q_{z*}} \;=\;
   \left\{ {0 \qquad\;\;\qmbox{if} \l(z)=m \qmbox{or} q_{z0}=\fr12 \atop
            1+N_{\v q_{z0*}}+N_{\v q_{z1*}} \;\;\qquad\qquad\qmbox{else}}
   \right.
\eeq
The effective dimension is zero if $q_{z0}=\fr12$, since this
implies that the whole tree $\G_z$ has $q_{zy}=\fr12$ due to
(\ref{EqPqRec}). If $q_z\neq\fr12$, we add the effective dimensions
of subtrees $\G_{z0}$ and $\G_{z1}$ to the root degree of
freedom $q_{z0}=1-q_{z1}$.
We may be interested in the expected effective dimension $E[N_{\v
q_*}|D]$. Inserting (\ref{eqNrec}) ($f(\v q_{z*})=N_{\v q_{z*}}$)
and (\ref{eqprec}) into the r.h.s.\ of (\ref{eqfexp}) we can
prove the following recursion for the expected effective dimension
\beq \label{eqEMD}
   E_z[N_{\v q_{z*}}\!|\!D_z]
   = g_z(\!D_z\!)[1\!+\!E_{z0}[N_{\v q_{z0*}}\!|\!D_{z0}]\!+\!E_{z1}[N_{\v q_{z1*}}\!|\!D_{z1}]]
\eeq
Read: The expected dimension of $\v q_{z*}$ (l.h.s.) equals to 1
for the root degree of freedom plus the expected dimensions in the
left and right subtrees, multiplied with the probability
$g_z(D_z)$ of splitting $\G_z$ (r.h.s.).
The recursion terminates with $E_z[N_{\v q_{z*}}|D_z]=0$ when
$\l(z)=m$.
Higher (central) moments like the variance can be computed
similarly. One can also compute the whole distribution $(P[N_{\v
q_*}=k|D])_{k\in\SetN}$ by convolution. Inserting
\beq\label{eqfdN}
  f(\v q_{z*})
  \;=\; \delta_{N_{\v q_{z*},k+1}}
  \;=\;
  \left\{ { 0 \!\qmbox{if}\! l=m \!\qmbox{or}\! q_{z0}=\fr12, \atop
  \sum_{i=0}^k \delta_{N_{\v q_{z0*},i}}\delta_{N_{\v q_{z1*},k-i}}
  \quad\mbox{else,} }
  \right.
\eeq
and (\ref{eqprec}) into (\ref{eqfexp}) we get
\bqa\nonumber
  & & \nq\nq P_z[N_{\v q_{z*}}=0|D_z] \;=\;
   1-g_z(D_z), \qquad \qmbox{for} l<m,
\\[1ex] \label{eqMDR}
  & & \nq\nq P_z[N_{\v q_{z*}}=k+1|D_z] \;=\;
\\[-1ex] \nonumber
  & & \nq\nq g_z(D_z)\!\cdot\!\!\sum_{i=0}^k
  P_{z0}[N_{\v q_{z0*}}\!=\!i|D_{z0}]\cdot P_{z1}[N_{\v q_{z1*}}\!=\!k\!-\!i|D_{z1}],
\\ \nonumber
  & & \nq\nq P_z[N_{\v q_{z*}}=k|D_z]=\delta_{k0}
  :=\textstyle\big\{{1 \;{\rm if}\; k=0 \atop 0 \;{\rm if}\; k>0}\big\} \qmbox{for} l=m.
\vspace{-2ex}
\eqa
Read: The probability that tree $\G_z$ has dimension $k+1$ equals
the probability of splitting, times the probability that left subtree
has dimension $i$, times the probability that right subtree has
dimension $k-i$, summed over all possible $i$. Again, this follows from
inserting (\ref{eqfdN}) and (\ref{eqprec}) into (\ref{eqfexp}).

Let us define a cell or bin as a maximal volume on which $q()$ is
constant. Then the model dimension is 1 less than the number of
bins (due to the probability constraint). Hence adding 1 to the
above quantities we also have expressions for the expected number
of cells and distribution.

\paradot{Tree height and cell size}
The effective height of tree $\v q_{z*}$ at $x\in\G_z$
is also an interesting property.
If $q_{z0}=\fr12$ or $\l(z)=m$, then the height $h_{\v q_{z*}}(x)$ of
tree $\v q_{z*}$ at $x$ is obviously zero. If
$q_{z0}\neq\fr12$, we take the height of the subtree
$\v q_{zx_{l+1}*}$ that contains $x$ and add 1:
\beqn\label{eqhrec}
   h_{\v q_{z*}}(x) \;=\;
   \left\{ {0 \atop 1+h_{\v q_{zx_{l+1}*}}(x)}
           {\mbox{if}\quad \l(z)=m \qmbox{or} q_{z0}=\fr12 \atop
             \mbox{else}\hspace{10ex}}   \right.
\eeqn
One can show that the tree height at $x\in\G_z$ averaged over all trees $\v q_{z*}$ is
\beqn
  E_z\![h_{\v q_{z*}}\!(x)|D_z] \;=\;
  g_z(D_z)\Big[1 + E_{zx_{l+1}}\![h_{\v q_{zx_{l+1}\!*\!}}(x)|D_{zx_{l+1}}]\Big]
\eeqn
We may also want to compute the tree height averaged over all
$x\in\G_z$.
\beqn
  \bar h_{\v q_{z*}}
  \;:=\;
  \int h_{\v q_{z*}}(x)q(x|\G_z)dx
  \;=\;
  \left\{ {0 \atop 1+q_{z0}\!\cdot\!\bar h_{\v q_{z0*}}+q_{z1}\!\cdot\!\bar h_{\v q_{z1*}}}
             {\mbox{if}\quad \l(z)=m \qmbox{or} q_{z0}=\fr12 \atop
               \mbox{else}\hspace{10ex}}   \right.
\eeqn\vspace{-2ex}
\bqan
  E_z[\bar h_{\v q_{z*}}|D_z] \;=\;
  g_z(D_z)\Big[1 \nq\;&+&\nq\; {n_{z0}\!+\a\over n_z+2\a}E_{z0}[\bar h_{\v q_{z0*}}|D_{z0}]
\\
            \nq\;&+&\nq\; {n_{z1}\!+\a\over n_z+2\a}E_{z1}[\bar h_{\v q_{z1*}}|D_{z1}]\Big]
\eqan
with obvious interpretation: The expected height of a subtree is
weighted by its relative importance, that is (an estimate of) its
probability. The recursion terminates with $E_z[h_{\v
q_{z*}}|D_z]=0$ when $\l(z)=m$. We can also compute intra and
inter tree height variances.

Finally consider the average cell size or volume $v$. Maybe more
useful is to consider the logarithm $-\lb|\G_z|=\l(z)$, since
otherwise small volumes can get swamped in the expectation by a
single large one. Log-volume $v_{\v q_{z*}}=\l(z)$ if $\l(z)=m$ or
$q_{z0}=\fr12$, and else recursively $v_{\v q_{z*}}=q_{z0}v_{\v
q_{z0*}}+q_{z1}v_{\v q_{z1*}}$. We can reduce this to the tree height,
since $v_{\v q_{z*}}=\bar h_{\v q_{z*}}+\l(z)$, in particular
$v_{\v q_*}=\bar h_{\v q_*}$

\paradot{\boldmath Moments in $x$}
Assume data $x\in\G=[0,1)$ are sampled from $q()$. The mean and
variance of $x$ w.r.t.\ $q()$ are important statistical
quantities. More generally consider
\beqn
  f(\v q_{z*})
  \;=\; M_{\v q_{z*}} \;:=\;
  {1\over|\G_z|}\int_{\G_z} M(x) q_z(x) dx
\eeqn
(${1\over|\G_z|}\int q_z(x) dx=1$).
Since $q()$ is itself random, it is natural to consider the $p$-expected
$q$-mean (\ref{eqfexp}) of $M$
\beq\label{eqEMom}
  \EE_z[M(x)|D_z]
  \;:=\; E_z[M_{\v q_{z*}}|D_z]
  \;=\; {1\over|\G_z|}\int_{\G_z} M(x) p_z(x|D_z) dx
\eeq
Inserting recursion (\ref{eqEvRec}) for $p_z(D_z,x)$ using
(\ref{eqPredDistr}) in (\ref{eqEMom}) we get
\beq\label{recEMom}
  \EE_z[M(x)|D_z] = {u\!\cdot\!\bar M_z\over p_z(D_z)} + g_z(D_z)
\Big[{n_{z0}\!+\a\over n_z+2\a}\EE_{z0}[M(x)|D_{z0}]
                   + {n_{z1}\!+\a\over n_z+2\a}\EE_{z1}[M(x)|D_{z1}]\Big],
\eeq
again with obvious interpretation: The expectation of $M$ is a
mixture of a uniform expectation $\bar
M_z:={1\over|\G_z|}\int_{\G_z}M(x)dx$ and the weighted average of
expectations in left and right subtree. The recursion terminates
with
$\EE_z[M(x)|D_z]=\bar M_z$ when $l=m$.
Examples: For $M(x)\equiv 1$, both sides of (\ref{recEMom})
evaluate to 1 as it should. For the $k$th moment of $x$,
$M(x)=x^k$ we have
$\bar M_z = [(z+1)^{k+1}-z^{k+1}]/[2^{kl}(k+1)]$,
where $z=2^l 0.z$ is interpreted as an integer in binary
representation. The distribution function $P_z[x\leq a|D_z]$ is
obtained for $M(x)=\{ {1\text{ if } x\leq a\atop 0 \text{ if }
x>a}$. For $a\in\G_z$ we have $\bar M_z=2^l a-z$. Since
$M(x)$ is constant on $\G_z\not\ni a$, we have
$\EE_z[M(x)|D_z]=\{{1 \text{ if } a\geq 0.z+2^{-l}\atop 0\text{ if
} a<0.z\phantom{+2^{-l}}}$ in this case, hence only one recursion
in (\ref{recEMom}) needs to be expanded (since $a\not\in\G_{z0}$
or $a\not\in\G_{z1}$).

\section{\boldmath Infinite Trees ($m\to\infty$)}\label{secIT}

\paradot{Motivation}
We have chosen an (arbitrary) finite tree height $m$ in our setup,
needed to have a well-defined recursion start at the leaves of the
trees. What we are really interested in are infinite trees
($m=\infty$).
Why not feel lucky with finite $m$? First, for continuous domain
$\G$ (e.g.\ interval $[0,1)$), our tree model contains only
piecewise constant models. The true distribution $\dot q()$ is
typically non-constant and continuous (Beta, normal, ...). Such
distributions are outside a finite tree model class (but inside
the infinite model), and the posterior $p(x|D)$ cannot converge to
the true distribution, since it is also piecewise constant. Hence
all other estimators based on the posterior are also not
consistent. Second, a finite $m$ violates scale invariance (a
non-informative prior on $\G_z$ should be the same for all $z$,
apart from scaling). Finally, having to choose the ``right'' $m$
may be worrisome.

For increasing $m$, the cells $\G_x$ become smaller and will
(normally) eventually contain either only a single data item, or be
empty. It should not matter whether we further subdivide empty or
singleton cells. So we expect inferences to be independent of $m$
for sufficiently large $m$, or at least the limit $m\to\infty$ to
exist. In this section we show that this is essentially true, but
with interesting exceptions.

\paradot{\boldmath Prior inferences ($D=\es$)}
We first consider the prior (zero data) case $D=\es$. Recall that
$z\in\SetB_0^m$ is some node and $x\in\SetB^m$ a leaf node.
Normalization implies $p_z(\es)=1$ for all $z$, which is independent
of $m$, hence the prior evidence exists for $m\to\infty$ (see below
for a formal proof). This is nice, but hardly surprising.

The prior effective model dimension $N_{\v q_*}$ is more
interesting. $D=\es$ implies $D_z=\es$ implies $n_z=0$ implies
$w(n_{z0},n_{z1})=1$ implies a chance $g_z(\es)=1-u=s$
for a split (see (\ref{eqSplitProb})).
The recursion (\ref{eqEMD}) reduces to
\beqn
   E_z[N_{\v q_{z*}}]
   = s[1+E_{z0}[N_{\v q_{z0*}}]+E_{z1}[N_{\v q_{z1*}}]]
\eeqn
with $E_z[N_{\v q_{z*}}]=0$ for $l\equiv \l(z)=m$ and is easily solved:
\beqn
   E_z[N_{\v q_{z*}}]
   \;=\; s(1+2s(1+...))
   \;=\; s{1-(2s)^{m-l}\over 1-2s}
   \;=\; \left\{
   \begin{array}{r@{\qmbox{if}}c}
     \stackrel{m\to\infty}\longrightarrow{s\over 1-2s} & s<\fr12 \\
     \fr12(m-l) \stackrel{lin.}\to\infty                & s=\fr12 \\
     \stackrel{exp.}\to\infty                          & s>\fr12 \\
   \end{array} \right.
\eeqn
This can be understood as follows. Consider trees truncated at
nodes with uniform distribution. Assume that there are $k(l)$
nodes at height $l$. With probability $u$ the node is a leaf, and
with probability $s$ it has two children. So the expected number
of nodes at height $l+1$ is $k(l+1)=2s\cdot k(l)$. So the number
of nodes exponentially increases/decreases with $l$ for
$s>\fr12/s<\fr12$, which results in a infinite/finite total number
of nodes (=dimension). The linear divergence for $s=\fr12$
looks alerting (overfitting), but we now show that
the distribution exists and an infinite expectation is actually a
good sign. Recursion (\ref{eqMDR}) reads
\beqn
  P_z[N_{\v q_{z*}}=k+1] =
  s\sum_{i=0}^k P_{z0}[N_{\v q_{z0*}}\!=\!i]\cdot P_{z1}[N_{\v q_{z1*}}\!=\!k\!-\!i]
\eeqn
with $P_z[N_{\v q_{z*}}=0]=u$ for $l<m$ and $P_z[N_{\v
q_{z*}}=k]=\delta_{k0}$ for $l=m$. So the recursion terminates in
recursion depth $\min\{k+1,m-l\}$. Hence $P_z[N_{\v q_{z*}}=k+1]$
is the same for all $m>l+k$, which implies that the limit
$m\to\infty$ exists. Furthermore, recursion and termination are
independent of $z$, hence also $a_k:=P_z[N_{\v q_{z*}}=k]$. So we
have to solve the recursion
\beq\label{eqMDar}
   a_{k+1} = s\sum_{i=0}^k a_i\!\cdot\! a_{k-i} \qmbox{with} a_0=u
\eeq
The first few coefficients can be bootstrapped by hand, e.g.\ for
$s=\fr12$ we get ($\v a=\fr12$,$1\over 8$,$1\over 16$,$5\over
128$,$7\over 256$,$21\over 1024$,$33\over 2048$,...).
A closed form  can
also be obtained: Inserting (\ref{eqMDar}) into $f(x):=\sum_{k=0}^\infty
a_k x^{k+1}$ we get $f(x)=ux+sf^2(x)$ with solution
$f(x)={1\over 2s}[1-\sqrt{1-4sux}]$, which has Taylor expansion coefficients
\beqn
   a_k = 2u(-4su)^k\left({1/2\atop k+1}\right)
   = {u(4su)^k\over(k\!+\!1)4^k}\left({2k\atop k}\right)
   \sim {u(4su)^k\over\sqrt{\pi}}\; k^{-3/2}
\eeqn
For $s\leq\fr12$, $(a_k)_{k\in\SetN_0}$ is a well-behaved properly
normalized probability measure ($\sum_k a_k=f(1)=1<\infty$). For
$s<\fr12$ it decreases exponentially in $k$, implying that all
moments exist and indicating strong bias towards simple models. For
$s=\fr12$, $a_k\sim k^{-3/2}$ decreases polynomially in $k$, too
slow for the expectation $E[N_{\v q_*}]=\sum_k k\cdot a_k=\infty$ to
exist, but this is exactly how a proper non-informative prior on
$\SetN$ should look like: as uniform as possible, i.e.\ slowly
decreasing. Further, $P[N_{\v q_*}<\infty]=\sum_k a_k=1$ shows that
the effective dimension is almost surely finite, i.e.\ infinite
(Polya) trees have probability zero for $s\leq\fr12$.
On the other hand, for $s>\fr12$, we have $P[N_{\v
q_*}=\infty]=1-f(1)=2-{1\over s}$, i.e.\ a non-zero probability
for infinite trees. The reason why $a_k$ also exponentially
decreases in this case is that as deeper a tree grows as less
likely it stays finite.
These results are consistent with the expected model dimension.
They indicate a proper behavior of our model for $s\leq\fr12$
and in particular for $s=\fr12$.

For the tree height we have
$E_z[h_{\v q_{z*}}(x)]=0$ if $l=m$ and otherwise
\bqan
  E_z[h_{\v q_{z*}}(x)] &=&
  s[1+E_{zx_{l+1}}[h_{\v q_{zx_{l+1}*}}(x)]]
\\
  &=& s+s^2+...+s^{m-l}
\\
  &=& \left\{
  \begin{array}{c@{\;\toinfty{m}\;}c@{\qmbox{if}}c}
    s\!\cdot\!{1-s^{m-l}\over 1-s} & {s\over 1-s} & s<1, \\
    m-l & \infty & s=1, \\
  \end{array} \right.
\eqan
i.e.\ the prior expected height is large/small if the splitting
probability $s$ is large/small. The same holds for the expected
average height $E_z[\bar h_{\v q_{z*}}]\to{s\over 1-s}$. This is the
first case where the result is not independent of $m$ for large
finite $m$, but it converges for $m\to\infty$, what is enough for
our purpose.
Note that a finite expected tree height even for $1>s\geq\fr12$ is
consistent with an infinite model dimension, since $h=\infty$ only
for a vanishing fraction of tree branches $x\in\G$, i.e.\ for a
set of measure zero.

The prior moments $M(x)$ are easy to compute: Since
$p_z(x)\equiv 1$, we get $\EE_z[M(x)]=\bar M_z$.

\paradot{\boldmath Multi-points $D=(x^1,...,x^1)$}
The next situation we analyze is multiple points
$D=(x^1,...,x^1)$, where all data points are identical. For
continuous spaces and non-singular priors, the probability of such
an event is zero, so this scenario does not seem particulary
interesting, but:
When computing posteriors, $x$ is not chosen randomly but
deliberately, so in computing $p(D,x)$, $x$ could be equal to
$x^1$ (although again only with probability zero). When computing
higher moments we need $p(D,x,...,x)$ and definitely
encounter multiple points. Also, multi-points help to understand
the case when $x$ or $x^i$ comes very close to $x^1$. Also the
true prior may be singular causing $x^1=x^2$ with non-zero
probability. Finally, the multi-point case includes $n=1$, which we
have to analyze in any case.

$D_z$ is either $\es$ or $D$. $D_z=\es$ has been studied in the
last \S, so we assume $D_z=D$. Either $D_{z0}=\es$ or
$D_{z1}=\es$. W.l.g.\ we assume the latter. Then $D_{z0}=D$,
$n_{z0}=n$, $n_{z1}=0$, which implies
$w(n_{z0},n_{z1})=2^{-n}{\Ga(n+2\a)\Ga(\a)\over\Ga(n+\a)\Ga(2\a)}$.
Defining $\bar w:=s/w$, recursion (\ref{eqEvDens})
reduces to
\bqa\nonumber
  & & \nq p_z(D) \;=\; u+\bar w\!\cdot\!p_{z0}(D) \;=...=\;
  u{1-\bar w^{m-l}\over 1-\bar w} + \bar w^{m-l}
\\ \nonumber
 & & \nq \left\{
 \begin{array}{lcl}
   \equiv 1                                       & \mbox{if} & \bar w=s \\
   \to {u\over 1-\bar w}                          & \mbox{if} & \bar w<1 \\
   = u(m\!-\!l)+1\stackrel{lin.}\longrightarrow\infty & \mbox{if} & \bar w=1 \\
   \to {\bar w-s\over\bar w-1}\bar w^{m-l}\stackrel{exp.}\longrightarrow\infty & \mbox{if} & \bar w>1 \\
 \end{array}
 \right.
\\ \label{eqRecMp}
\eqa
So the evidence exists for $m\to\infty$ iff $\bar w<1$.
For $n=0$ and $n=1$ we have $\bar w=s<1$ (excluding $s=1$), hence
$p(\es)\equiv 1$ is correctly normalized, as claimed in the
previous \S, and $p(x)\equiv 1$ is
uniform as symmetry demands.
For double points $n=2$ the evidence $p(x,x)$ is still finite iff
$\bar w=s\cdot{\a+1\over\a+1/2}<1$. The latter is true for
$s\leq\fr12$ and all $\a>0$.
For any $n$ but $\a=1$, $\bar w<1$ iff $s<(n+1)2^{-n}$, i.e.\
higher multi-point evidences only exist for exponentially small
$s$. If we want $\bar w<1$ $\forall\a>0$, $s<2^{1-n}$ has to
be even smaller. This follows from $\bar w$ being increasing in
$n$ and decreasing in $\a$ ($\bar w\searrow s$ for $\a\to\infty$).
To conclude: For every fixed $\a$ and $s$, multiplicity $n$ of
points must not be too high, but for any $n$ one can choose
$\a$ sufficiently large or $s$ sufficiently small so that the
$n$-multi-point evidence and hence the $n^{th}$ moments
exist.

For the mean of $M(x)$ we get
\beqn\textstyle
  \EE_z[M(x)|D] \;=\; (1\!-\!\tilde w)\bar M_z + \tilde w{\a\over n+2\a}\bar M_{z0}
  + (\tilde w{n+\a\over n+2\a})\EE_{z1}[M(x)|D]
\eeqn
if $x^1\in\G_{z1}$ and similarly for $x^1\in\G_{z0}$, where
$\tilde w=\min\{\bar w,1\}$. This is not a self-consistency
equation, but since $\tilde w{n+\a\over n+2\a}<1$, the linear
recursion converges exponentially to the exact value.

\paradot{\boldmath Multi-points for $\a=1$ and $s=\fr12$}
We present some more results for the most interesting case $\a=1$
and $s=\fr12$, which we will also investigate numerically.

We see that $p(\es)\equiv 1$ is correctly normalized, and
$p(x)\equiv 1$ is uniform as symmetry demands. For double points,
the evidence $p(x,x)\to{3\over 2}$ is still finite. It diverges
linearly for triple points and exponentially for quadruple-and-higher
points. So $q(x)$ has finite prior mean $E[q(x)]=1$ and
variance $\mbox{Var}[q(x)]=\fr12$ (Section \ref{secTMM}). The
skewness and kurtosis are infinite, indicating a heavy tail, as
desired for a non-informative prior.

Since $p(D)\equiv 1$ and $w=1$ for $n=1$ are the same as for the
$n=0$ case, all prior $n=0$, $m\to\infty$ results remain valid for
$n=1$: $g(x)=\fr12$, $E[N_{\v q_*}|x]=\infty$, $P[N_{\v
q_*}=k|x]=a_k$, and $E[h_{\v q_*}(x)|x]\to 1$.

For $n=2$ we get $g(x,x)={2\over 3}$, $b_k:=P[N_{\v q_*}=k|x,x]$,
$b_0=1-g(x,x)={1\over 3}$, $b_{k+1} = {2\over 3}\sum_{i=0}^k
b_i\!\cdot\! a_{k-i} =$ ($1\over 9$,$7\over 108$,$29\over
648$,...), $h(x):=\sum_{k=0}^\infty b_k x^{k+1} = {1\over
3}[1+2h(x)f(x)]={x\over 1+2\sqrt{1-x}}$, and $E[h_{\v
q_*}(x)|x,x]\to 2$.

For $n\geq 3$ we have $g_z(D)=1$, $b_0=0$ $\Rightarrow$
$b_k\equiv 0$ $\Rightarrow$ $P[N_{\v q_*}<\infty|D]=0$, $E[h_{\v
q_*}(x)|D]=m\to\infty$ for $x=x^1$. The tree at $x$ has
infinite height and singular distribution.

For all $n$ we have $p(D)\geq 1$, $g(D)\geq\fr12$, $E[N_{\v q_*}|D]\sim 2
g(D)$ for $m\to\infty$,

Not a double point, but also straightforward to compute is
$p(x,y)={3\over 2}-({2\over 3})^{l+1}$ if $x\in\G_{z0}$ and
$y\in\G_{z1}$, i.e.\ $x$ and $y$ separate at level $l=\l(z)$,
consistent with $p(x,x)={3\over 2}$ ($l\widehat=\infty$).

\begin{figure}
\centerline{\includegraphics[width=0.5\textwidth]{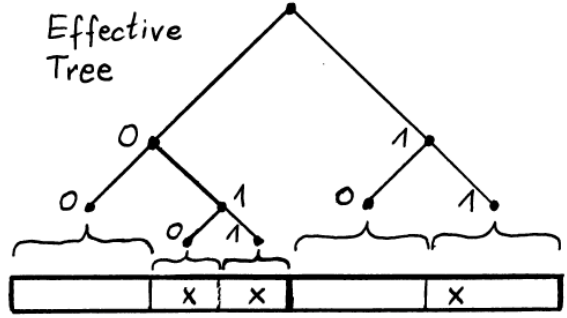}}
\caption{\label{figEffectiveTree}Effectively recursed tree.
Closed form solutions are used for intervals containing no or only
a single (multi)point.}
\end{figure}

\paradot{\boldmath General $D$}
We now consider general $D$. In order to compute $p(D)$ and other
quantities, we recurse (\ref{eqEvDens}) down the tree until $D_z$
is a multi-point $D_z=D_z^= :=(x',...,x')$ with $x'\in\G_z$. We
call the depth $m_{x'}:=\l(z)$ at which this happens, the
separation level. If we consider $n_z\in\SetN_0$, this also
includes the most important empty and singleton case. In this way,
the recursion always terminates. For instance, for $\G=[0,1)$, if
$\eps:=\min\{|x^i-x^j|:x^i\neq x^j$ with $x^i,x^j\in D\}$ is the
shortest distance, then $m_{x'}<\lb{2\over\eps}=:m_0<\infty$,
since $\eps>0$. At the separation level we can insert the derived
formulas (for evidence, posterior, dimension, height, moments) for
multi-points (Figure \ref{figEffectiveTree}). Note, there is no
approximation here. The procedure is exact, since we analytically
computed the infinite recursion for multi-points.

So we have devised a finite procedure for exactly computing all
quantities of interest. In the worst case, we have to recurse down
to level $m_0$ for each data point, hence our procedure has
computational complexity $O(n\cdot m_0)$. For non-singular prior,
the time is actually $O(n)$ with probability 1. So, inference in
our mixture tree model is {\em very} fast. Polya trees for
suitable Beta prior (should) admit similar algorithms.

\paradot{Multi-point divergences}
Consider again $s=\fr12$ and $\a=1$. Since $E_z[N_{\v
q_{z*}}|D_z^=]=\infty$ recurses up, we have $E[N_{\v
q_{*}}|D]=\infty$ for all $D$. We now discuss possible divergences
caused by true multi-points $n_z^=>2$ in $D$. $E_z[h_{\v
q_{z*}}(x)|D_z^=]=\infty$ recurses up to $E[h_{\v
q_{*}}(x)|D]=\infty$. Similarly, $P_z[N_{\v
q_{z*}}<\infty|D_z^=]=0$, recurses up to $P[N_{\v
q_{*}}<\infty|D]=0$, since $g_z(D_z)=1$ along the path.

There are interesting cases where $p(D)=\infty$, but posteriors
are finite, since infinities cancel out. These are $p(x|D)$ if
$x$ occurs at most once in $D$, and $p(x,x,|D)$ if $x\not\in D$, otherwise
$p(x,...,x|D)=\infty$. This is is very welcome: $E[q(x)|D]<\infty$ for
all $x$, if $D$ contains only singletons, which is true w.p.1 for
all regular $\dot q()$. The posterior variance of $q(x)$ is
finite iff $x\not\in D$, which could be better.
We adapt recursion (\ref{eqEvDens}) by scaling $p(D)$, $p(D,x)$
and $p(D,x,x)$ with the same constant $c_m\to 0$ such that they
stay finite, which works since $D$, $(D,x)$ and $(D,x,x)$ have the
same triple, quadruple, ... points. Choose $l$ large enough
(separation level) so that $D_z$ is a triple-or-higher point, and
w.l.g.\ $D_{z1}$ is empty. Then the recursion (\ref{eqEvDens})
reduces to recursion (\ref{eqRecMp}). The $\fr12+$ in
(\ref{eqEvDens}) gets swamped by $p_{z0}(D_{z0}^=)=\infty$ and can
be dropped. The remaining recursion is just a multiplication with
$\bar w$, which allows us to rescale $p_z(D_z^=)$ to
$p_z(D_z^=)=\bar w^{-l}$ (cf.\ (\ref{eqRecMp})). We return a flag in
the recursion indicating that the $\fr12+$ shall be dropped along
the path back to the root, since the true original $p_z(D_z^=)$
was infinite and would have swamped them all. We compute $p(D)$,
$p(D,x)$, and $p(D,x,x)$ in this way to arrive at finite
posteriors (not forgetting that the numbers we get for the
evidences themselves are fictitious).

A heuristic way of regularizing our model to yield always finite
results could be to cut the recursion short at the separation level
by definition, and then assign some regular (e.g.\ uniform)
distribution to this leaf. Since the most important quantities are
finite anyway, we refrain from such a data-dependent non-Bayesian hack.
Better is to assign a smaller prior weight $s<\fr12$ to a split;
then more (higher) moments become finite.

\paradot{\boldmath Consistency ($n\to\infty$)}
What remains to be shown is posterior consistency for $m=\infty$
similarly to the $m<\infty$ case. We will show weak consistency in
the sense that $p(\{q_z\}|D)$ concentrates around $\{\dot q_z\}$,
where $\{q_z\}$ is a {\em finite} collection of branching
probabilities. Consistency holds because the recursion for
$p(\{q_z\}|D)$ terminates at a depth independent of $m$ (for
sufficiently large $m$), so we are effectively in the finite tree
case. The only difficulty is that the recursion involves
$p_z(D_z)$ which still has recursion depth (i.e.\ depends on) $m$.
One solution could be to assume that all observations have some
finite precision $2^{-m'}$ and $x\in D\leadsto \G_{x'}\in D'$,
where $\l(x')=m'$ and $x\in\G_{x'}$. This would make all involved
recursions terminate at depth $m'$ and hence all recursions and
results for finite $m$ apply (with $m'$ instead of $m$).
More interesting is to keep $D$ and treat the $m=\infty$ case
properly. We show that for $n\to\infty$, the evidence and the
posterior converge uniformly in $m<\infty$, which implies
convergence also for $m=\infty$.
We sometimes indicate the $m$-dependence of $p_z$ by $p_z^m$ and
define $p_z^\infty:=\lim_{m\to\infty} p_z^m$ if the limit exists
(possibly infinite), but mostly drop the superscript
$m\leq\infty$.

\begin{theorem}[Weak consistency for infinite trees]\label{thWCIT}
Let $\l(z)=l\leq m'<m+1\leq\infty$ and $\v q_{z*'} =
(q_{zy}:y\in\SetB_1^{m'-l})$. Then evidence and p(oste)rior
$p_z(...)$ exist for $m=\infty$ and have the following properties
for $m\leq\infty$, where convergence $\toinfty{n_z}$ holds w.p.1
and is uniform in $m$:
\begin{list}{}{}
\item[i)] The marginal prior $p_z(\v q_{z*'})$ is independent of $m$,
\item[ii)] $p_z(D_z)\;\mathop{\longrightarrow}\limits^{n_z\to\infty}_{exp.}\;\infty$
for $us>0$ and non-uniform $\dot q_z()$.
\item[iii)] $\displaystyle
  p_z(\v q_{z*'}|D_z) \quad
  \mathop{\longrightarrow}\limits^{n_z\to\infty}_{weak}
  \prod_{y\in\SetB_1^{m'-l}} \delta(q_{zy}-\dot q_{zy})$.
\item[iv)] $p_z(\G_{\!y^1},...,\G_{\!y^k}|D_z) \;\toinfty{n_z}\;
\dot q_z(\G_{\!y^1})\cdot...\cdot\dot q_z(\G_{\!y^k})\quad$
for $y^i\in\SetB^*$ and $k\in\SetN$.
\end{list}
\end{theorem}
$(iv)$ implies weak convergence of $p(x|D)$ to $\dot q(x)$ in the
sense that $\int f(x)p(x|D)dx\to\int f(x)\dot q(x)dx$ for
continuous functions $f$, by an argument similar to the proof of
\cite[Thm.2.2]{Fabius:64}. Also the distribution function $P[x\leq
a|D]\to P[x\leq a|\dot q]$.

\paradot{Proof}
We only have to consider $m<\infty$. The $m=\infty$ case follows for
$(i)$ by definition and for $(ii)-(iv)$ since convergence
$\toinfty{n_z}$ is uniform in $m$.

{\boldmath $(i)$} The prior marginal is
\beq\label{eqPrMar}
  p_z(\v q_{z*'}) \;=\;  \int p_z(\v q_{z*})
  \nq\nq\prod_{y\in\SetB_1^{m-l}\setminus\SetB_1^{m'-l}}\nq\nq dq_{zy}
\eeq
For $l=m'$ we have $p_z(\v q_{z*'})=p_z(\es)=1$ independent of
$m$. For $l<m'$, inserting recursion (\ref{EqPqRec}) into (\ref{eqPrMar}) and
using $\SetB_1^{m-l}\setminus\SetB_1^{m'-l} =
\{0,1\}\times(\SetB_1^{m-(l+1)}\setminus\SetB_1^{m'-(l+1)})$ and
(\ref{eqPrMar}) for $z0$ and $z1$ backwards, we get
\beqn
  p_z(\v q_{z*'}) = u\!\cdot\!\nq\prod_{y\in\SetB_1^{m'-l}}\nq\; \delta(q_{zy}\!-\!\fr12)
  + s\!\cdot\!\Beta(q_{z0},q_{z1}|\a,\a)p_{z0}(\v q_{z0*'})p_{z1}(\v q_{z1*'})
\eeqn
So the recursion of $p_z^m(\v q_{z*'})$ and its termination is
independent of $m$, hence $p_z^m(\v q_{z*'})$ is independent of
$m$, hence $p_z^\infty(\v q_{z*'})$ exists and equals $p_z^m(\v
q_{z*'})$ for $m\geq m'$.

{\boldmath $(ii)$} First note that $p_z(D_z)$ in general depends on $m$.
Further, $p_z(D_z)\geq u>0$ $\forall z$. Assume first that $\dot q_{z0}\neq\fr12$.
Then from (\ref{eqEvRec}) and Theorem \ref{wAsym}$(iii)$ we get
\beq\label{eqDDivi}
    p_z(D_z) \;=\; u + s\!\cdot\!{p_{z0}(D_{z0})p_{z1}(D_{z1})\over w_{n_z}(\Delta_z)}
    \;\geq\; {su^2\over w_{n_z}}
    \;\mathop{\longrightarrow}\limits^{n_z\to\infty}_{exp.}\; \infty
\eeq
Divergence of $p_z(D_z)$ is uniform in $m$, since $w_{n_z}$ is
independent $m$. Now consider the more general case of non-uniform
$\dot q_z()$, i.e.\ $\exists y:\dot q_{zy0}\neq\fr12$. Then
(\ref{eqDDivi}) implies
$p_{zy}(D_{zy})\stackrel{exp.}\longrightarrow\infty$. Further, for
any $z$, if $p_{z0}\stackrel{exp.}\longrightarrow\infty$, then
$p_z\stackrel{exp.}\longrightarrow\infty$, since $p_{z1}\geq u>0$
and $w_{n_z}=O(\sqrt{n_z})$ by Theorem \ref{wAsym}, and similarly if
$p_{z1}\to\infty$. So by induction,
$p_{zy}\stackrel{exp.}\longrightarrow\infty$ uniformly in $m$
implies $p_z\stackrel{exp.}\longrightarrow\infty$ uniformly in $m$.
Finally, $p_z^\infty(D_z)$ exists, since $p_z^m(D_z)$ is independent
$m$ beyond the data separation level, as shown earlier.

{\boldmath $(iii)$} Similarly to the marginal prior recursion in
$(i)$ one can show that the recursion for the marginal posterior
$p_z(\v q_{z*'}|D_z)$ has the same form as the recursion
(\ref{eqprec}) of $p_z(\v q_{z*}|D_z)$ with $m$ replaced by $m'$.
Contrary to the prior, the posterior still depends on $m$ through
$p_z(D_z)$. Choosing $m$ beyond the data separation level does not
help since it increases with $n_z$. Nevertheless, the proof of
$(iii)$ is the same as for Theorem \ref{thpc} with $m\leadsto m'$.
Convergence is uniform in $m$, since convergence (\ref{eqpcr}) and
divergence of $p_z(D_z)$ are uniform. Finally, $p_z^\infty(\v
q_{z*'}|D_z)$ exists, since its recursion is finite (terminates at
$m'$) and $p_z^\infty(D_z)$ exists.

{\boldmath $(iv)$} Choose $m'\geq\max\{\l(y^1),...,\l(y^k)\}$. Then
\beqn
   p_z(\G_{\!y^1},...,\G_{\!y^k}|D_z)
   \;=\; \int q_z(\G_{\!y^1})\cdot...\cdot q_z(\G_{\!y^k})\cdot p_z(\v q_{z*'}|D_z) d\v q_{z*'}
\eeqn
exists, and $(iv)$ now follows from $(iii)$ and
$q_z(\G_y)=q_{y_1}\cdot...\cdot q_{y_{1:\l(y)}}$.
\qed

\section{The Algorithm}\label{secAlg}

\paradot{What it computes}
In the last two sections we derived all necessary formulas for
making inferences with our tree model. Collecting pieces together
we get the exact algorithm for infinite tree mixtures \ifalgin
below\else presented in Table \ref{BTAlg}\fi. It computes the
evidence $p(D)$, the expected tree height $E[h_{\v q_*}(x)|D]$ at
$x$, the average expected tree height $E[\bar h_{\v q_*}|D]$, and
the model dimension distribution $P[N_{\v q_*}|D]$. It also
returns the number of recursive function calls, i.e.\ the size of
the explicitly generated tree. The size is proportional to $n$ for regular
distributions $\dot q$.

\ifalgin
\paranodot{The BayesTree algorithm} (in pseudo C code)
\else
\begin{table}[htb]
\caption{\label{BTAlg}
{\bf BayesTree algorithm in pseudo C code}
\fi
takes arguments $(D[],n,x,N)$; data array $D[0..n-1]\in[0,1)^n$, a
point $x\in\SetR$, and an integer $N$. It returns $(p,h,\bar
h,\tilde p[],r)$; the logarithmic data evidence $p\widehat=\ln
p(D)$, the expected tree height $h\widehat=E[h_{\v q_*}(x)|D]$ at
$x$, the average expected tree height $\bar h\widehat=E[\bar h_{\v
q_*}|D]$, the model dimension distribution $\tilde
p[0..N-1]\widehat=P[N_{\v q_*}=..|D]$, and the number of recursive
function calls $r$ i.e.\ the size of the generated tree. $s$, $u$ and $\a$
are the global model parameters. Computation time is about $N^2 n\log
n$ nano-seconds on a 1GHz P4 laptop.
\ifalgin\else}\fi

\begin{list}{}{\parskip=0ex\parsep=0ex\itemsep=0.5ex\leftmargin=0ex\labelwidth=0ex}
  \item {\bf\boldmath BayesTree($D[],n,x,N$)}
  \begin{list}{}{\parskip=0ex\parsep=0ex\itemsep=0.5ex\leftmargin=2ex\labelwidth=1ex\labelsep=1ex}
    \item[$\lceil$] if ($n\leq 1$ and ($n==0$ or $D[0]==x$ or $x\not\in[0,1)$))
    \begin{list}{}{\parskip=0ex\parsep=0ex\itemsep=0.5ex\leftmargin=2ex\labelwidth=1ex\labelsep=1ex}
      \item[$\lceil$] if ($x\in[0,1)$) $h=s/u$; else $h=0$;
      \item     $\bar h=s/u$; $p=\ln(1)$; $r=1$;
      \item[$\lfloor$] for$(k=0,..,N-1)$ $\tilde p[k]=a_k$; \hfill /* see (\ref{eqMDar}) */
    \end{list}
      \item else
    \begin{list}{}{\parskip=0ex\parsep=0ex\itemsep=0.5ex\leftmargin=2ex\labelwidth=1ex\labelsep=1ex}
      \item[$\lceil$] $n_0=n_1=0$;
      \item           for$(i=0,..,n-1)$
      \begin{list}{}{\parskip=0ex\parsep=0ex\itemsep=0.5ex\leftmargin=2ex\labelwidth=1ex\labelsep=1ex}
        \item[$\lceil$] if ($D[i]<\fr12$) then[$\,D_0[n_0]=2D[i]$; \ \ \ $\,n_0=n_0+1$;]
        \item[$\lfloor$] \hspace{12ex} else [$D_1[n_1]=2D[i]-1$; $n_1=n_1+1$;]
      \end{list}
      \item ($p_0,h_0,\bar h_0,\tilde p_0[],r_0$)=BayesTree($D_0[],n_0,2x,N-1$);
      \item ($p_1,h_1,\bar h_1,\tilde p_1[],r_1$)=BayesTree($D_1[],n_1,2x\!-\!1,N\!-\!1$);
      \item $t=p_0+p_1-\ln w(n_0,n_1)$; \hfill /* see (\ref{eqWeights}) */
      \item if ($t<100$) then $p=\ln(u+s\cdot\exp(t))$;
      \item \hspace{10ex} \ else \ $p=t+\ln(s)$;
      \item $g=1-u\cdot\exp(-p)$;
      \item if ($x\in[0,1)$) then $h=g\cdot(1+h_0+h_1)$; else $h=0$;
      \item $\bar h=g\cdot(1+{n_0+\a\over n+2\a}\bar h_0+{n_1+\a\over n+2\a}\bar h_1)$;
      \item $\tilde p[0]=1-g$;
      \item for($k=0,..,N-1$) $\tilde p[k+1]=g\cdot\!\sum_{i=0}^k \tilde p_0[i]\cdot \tilde p_1[k-i]$;
      \item[$\lfloor$] $r=1+r_0+r_1$;
    \end{list}
    \item [$\lfloor$] {\bf\boldmath return ($p,h,\bar h,\tilde p[],r$); }
  \end{list}
\end{list}
\ifalgin\else\end{table}\fi

\paradot{How algorithm BayesTree() works}
Since evidence $p(D)$ and weight $1/w_n$ can grow exponentially
with $n$, we have to store and use their logarithms.
So the algorithm returns $p\widehat=\ln p(D)$.
In the $n\leq 1$ branch, the closed form solutions $p\widehat=\ln
p(\es)=\ln(1)$, $h\widehat=E[h_{\v q_*}(x)|\es\mbox{ or }x]=1$,
$\bar h\widehat=E[\bar h_{\v q_*}|D]=1$, and $\tilde
p[k]=a_k$ have been used to truncate the recursion. If $D=(x^1)\neq
x$, we have to recurse further until $x$ falls in an empty
interval.
In this case or if $n>1$ we partition $D$ into points left and
right of $\fr12$. Then we rescale the points to [0,1) and store
them in $D_0$ and $D_1$, respectively. Array $D$ could have been
reused (like in quick sort) without allocating two new arrays.
Then, algorithm BayesTree() is recursively called for each
partition. The results are combined according to the recursions
derived in Section \ref{secTMM}.
$\ln w$ can be computed from (\ref{eqWeights}) via $\ln
n!=\sum_{k=1}^n\ln k$. (Practically, pre-tabulating $a_k$ or $n!$
does not improve overall performance).
For computing $p$ we need to use $\ln(\fr12(1+e^t))\dot=t-\ln 2$ to
machine precision for large $t$ in order to avoid numerical overflow.

\paradot{Remarks}
Strictly speaking, the algorithm has runtime $O(n\log n)$, since
the sorting effectively runs once through all data at each level.
If we assume that the data are presorted or the counts $n_z$ are given,
then the algorithm is $O(n)$ \cite{Hutter:07btcode}.

We have not presented the part handling multi-points.
Given the formulas in Section \ref{secIT}, this is easy.
The complete C code is available from \cite{Hutter:07btcode}.

Note that $x$ passed to BayesTree() is {\em not} and cannot be
used to compute $p(x|D)$. For this, one has to call BayesTree()
twice, with $D$ and $(D,x)$, respectively. \ifalgin\else Computation time is
about $N^2 n\log n$ nano-seconds on a 1GHz P4 laptop.\fi The
quadratic order in $N$ is due to the convolution, which could be
reduced to $O(N\log N)$ by transforming it to a scalar product in
Fourier space with FFT.

Multiply calling BayesTree(), e.g.\ for computing the predictive
density function $p(x|D)$ on a fine $x$-grid, is inefficient. But
it is easy to see that if we once pre-compute the evidence
$p_z(D_z)$ for all $z$ up to the separation level in time $O(n)$,
we can compute ``local'' quantities like $p(x|D)$ at $x$ in time
$O(\log n)$. This is because only the branch containing $x$ needs
to be recursed, the other branch is immediately available, since it
involves the already pre-computed evidence only. The predictive
density $p(x|D)=E[q(x)|D]$ and higher moments, the distribution
function $P[x\leq a|D]$, updating $D$ by adding or removing one
data item, and most other local quantities can be computed in time
$O(\log n)$ by such a linear recursion.

A good way of checking correctness of the implementation {\em and}
of the derived formulas, is to force some {\em minimal} recursion
depth $m'$. The results must be independent of $m'$, since the
closed-form speedups are exact and applicable anywhere beyond the
separation level.

\section{Numerical Examples}\label{secEx}

\def\fourgraphs#1#2#3#4#5#6{
\begin{figure*}\centerline{
\includegraphics[width=0.5\columnwidth, height=0.16\textheight]{#1}
\includegraphics[width=0.5\columnwidth, height=0.16\textheight]{#2}
\includegraphics[width=0.5\columnwidth, height=0.16\textheight]{#3}
\includegraphics[width=0.5\columnwidth, height=0.16\textheight]{#4}}%
\vspace*{-2ex}\caption{\label{#5}#6}\vspace*{-4ex}
\end{figure*}}

\def\fourgraphs#1#2#3#4#5#6{
\begin{figure}\def\lfh{0.2}
\centerline{\includegraphics[width=0.5\textwidth, height=\lfh\textheight]{#1}
            \includegraphics[width=0.5\textwidth, height=\lfh\textheight]{#2}}
\centerline{\includegraphics[width=0.5\textwidth, height=\lfh\textheight]{#3}
            \includegraphics[width=0.5\textwidth, height=\lfh\textheight]{#4}}%
\caption{\label{#5}#6}
\end{figure}}

\fourgraphs{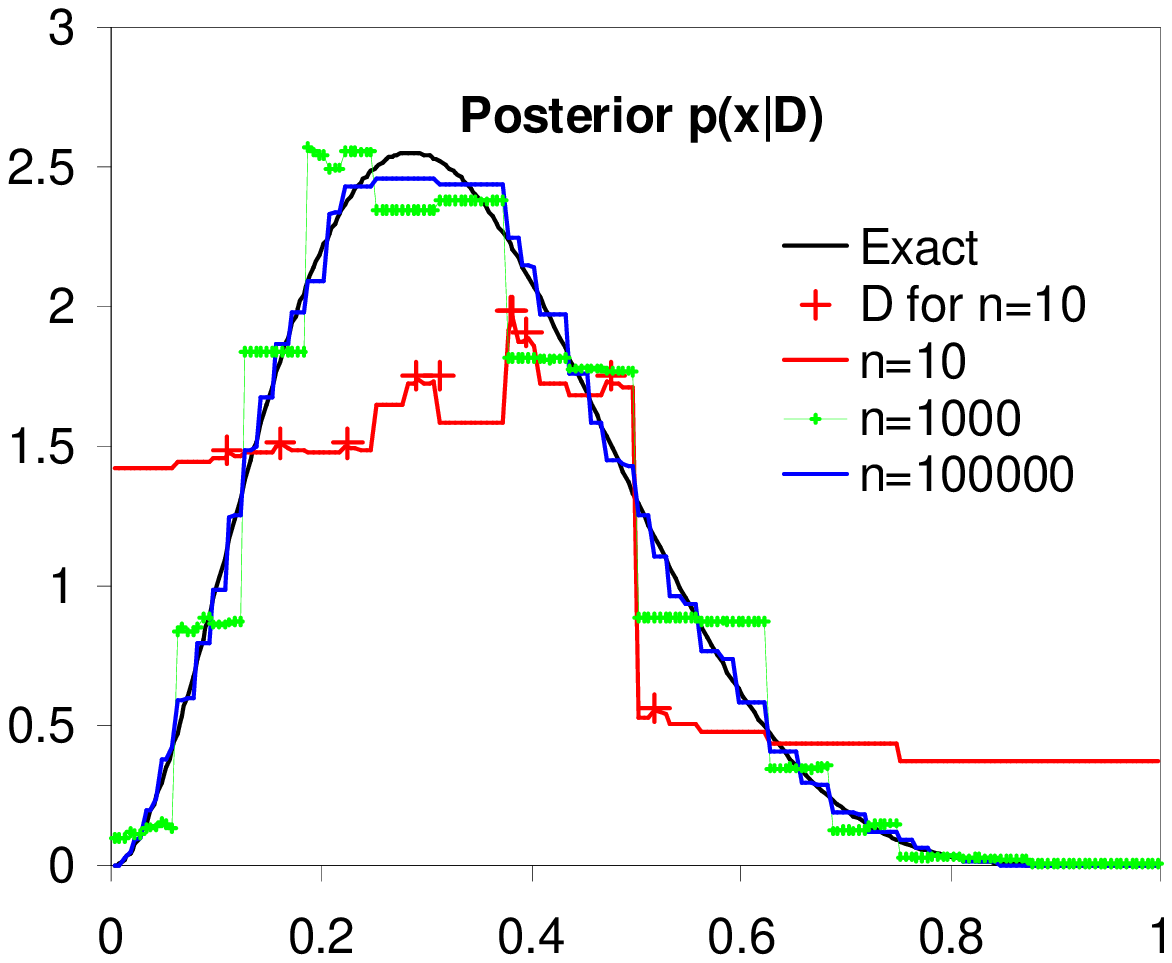}{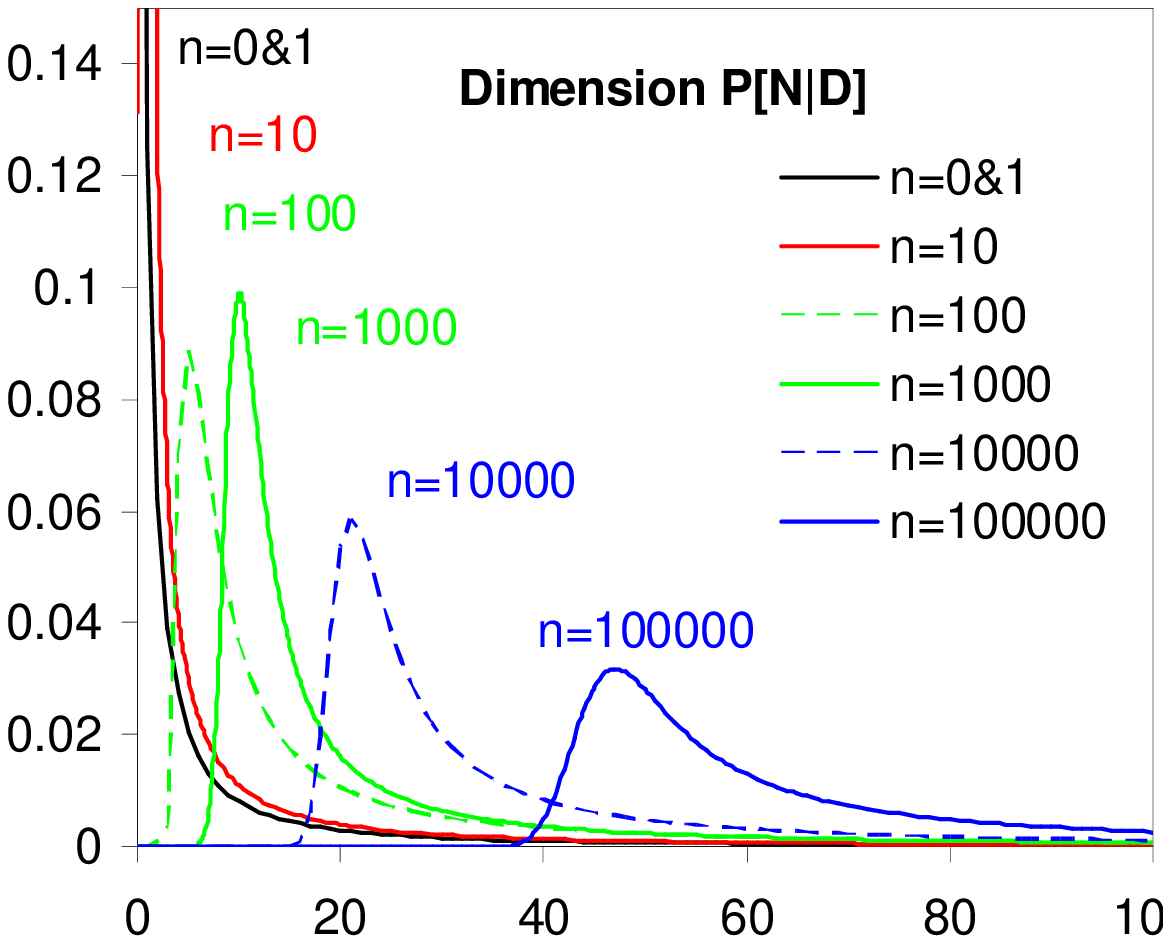}%
{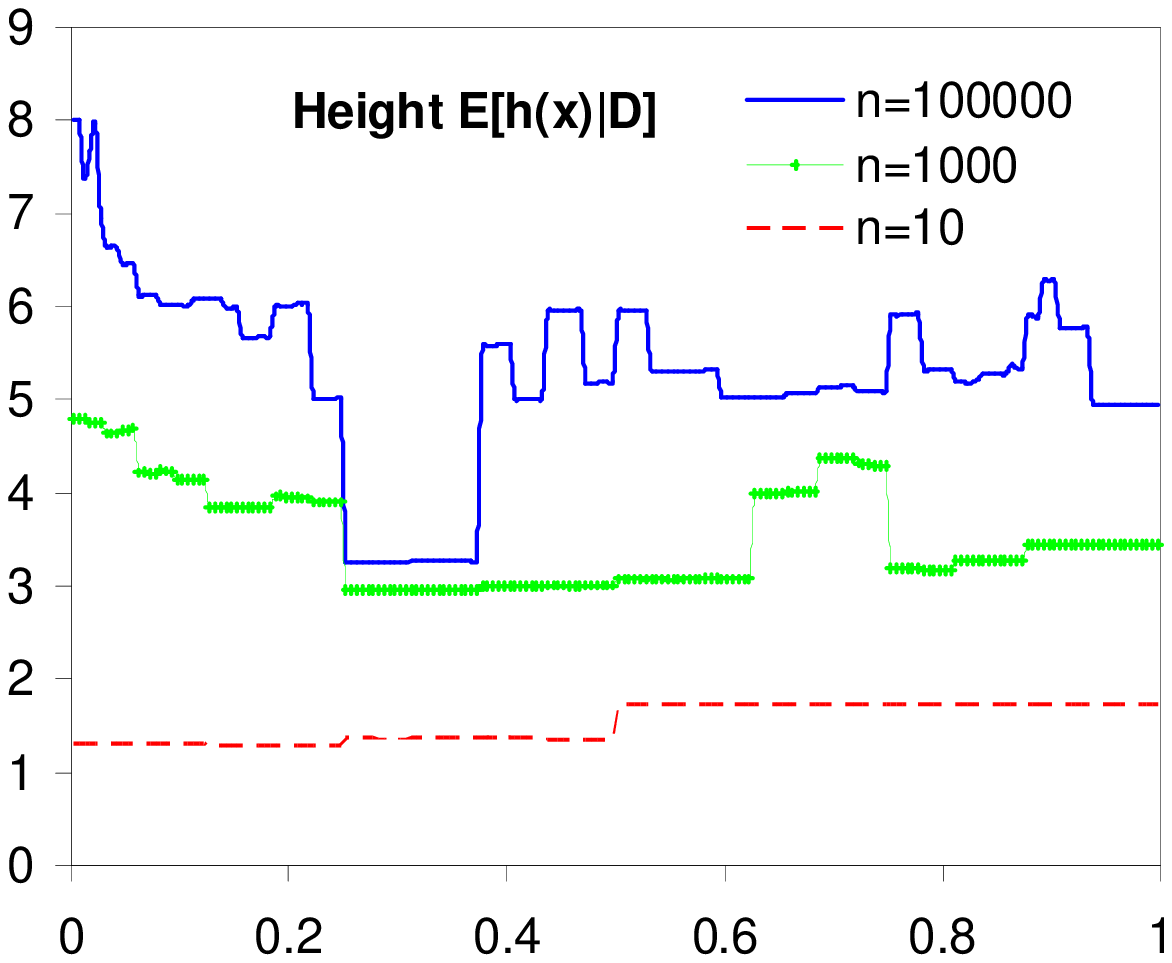}{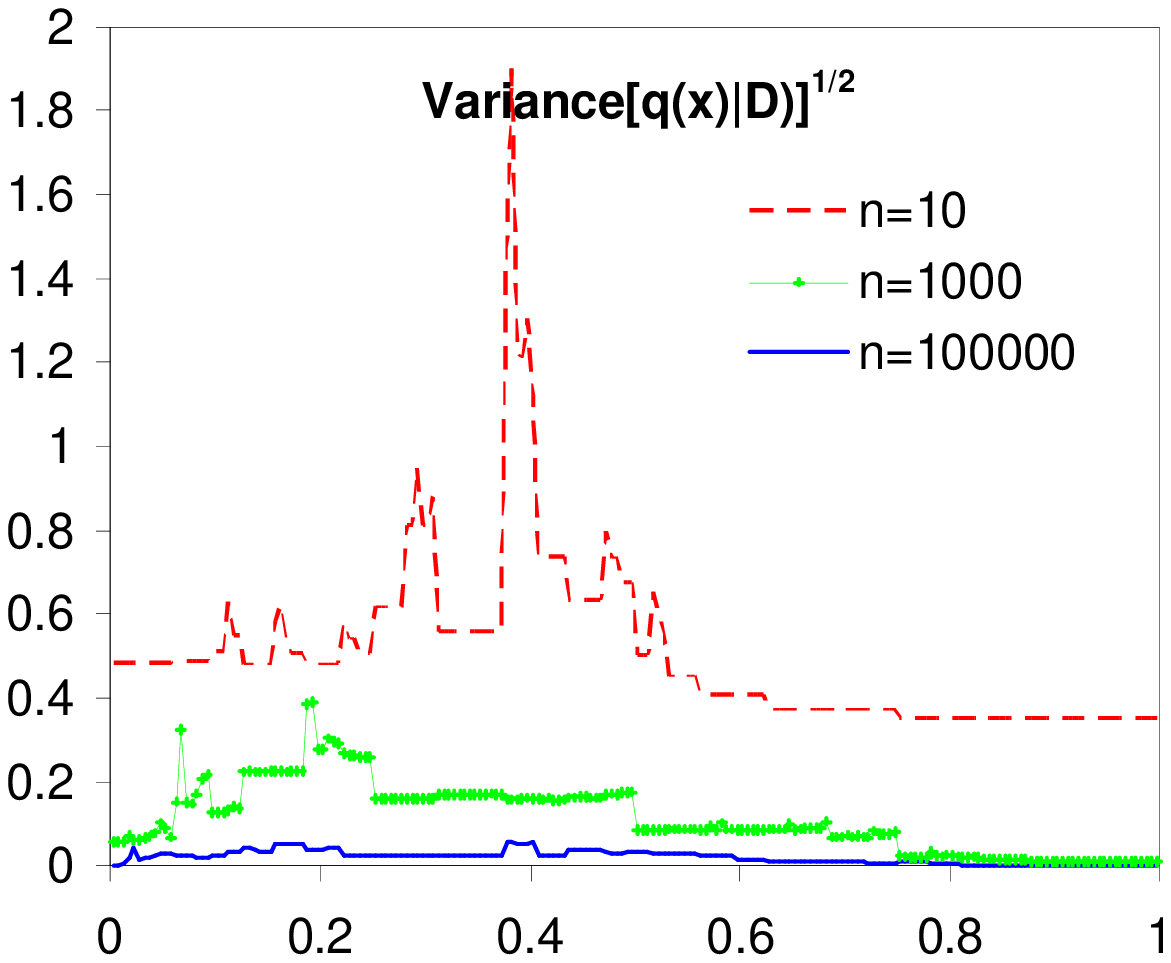}{figBeta}%
{BayesTree() results for the {\bf Beta(3,6)}$\propto {x^2(1-x)^5}$
distribution, prototype for a smooth distribution.}

\fourgraphs{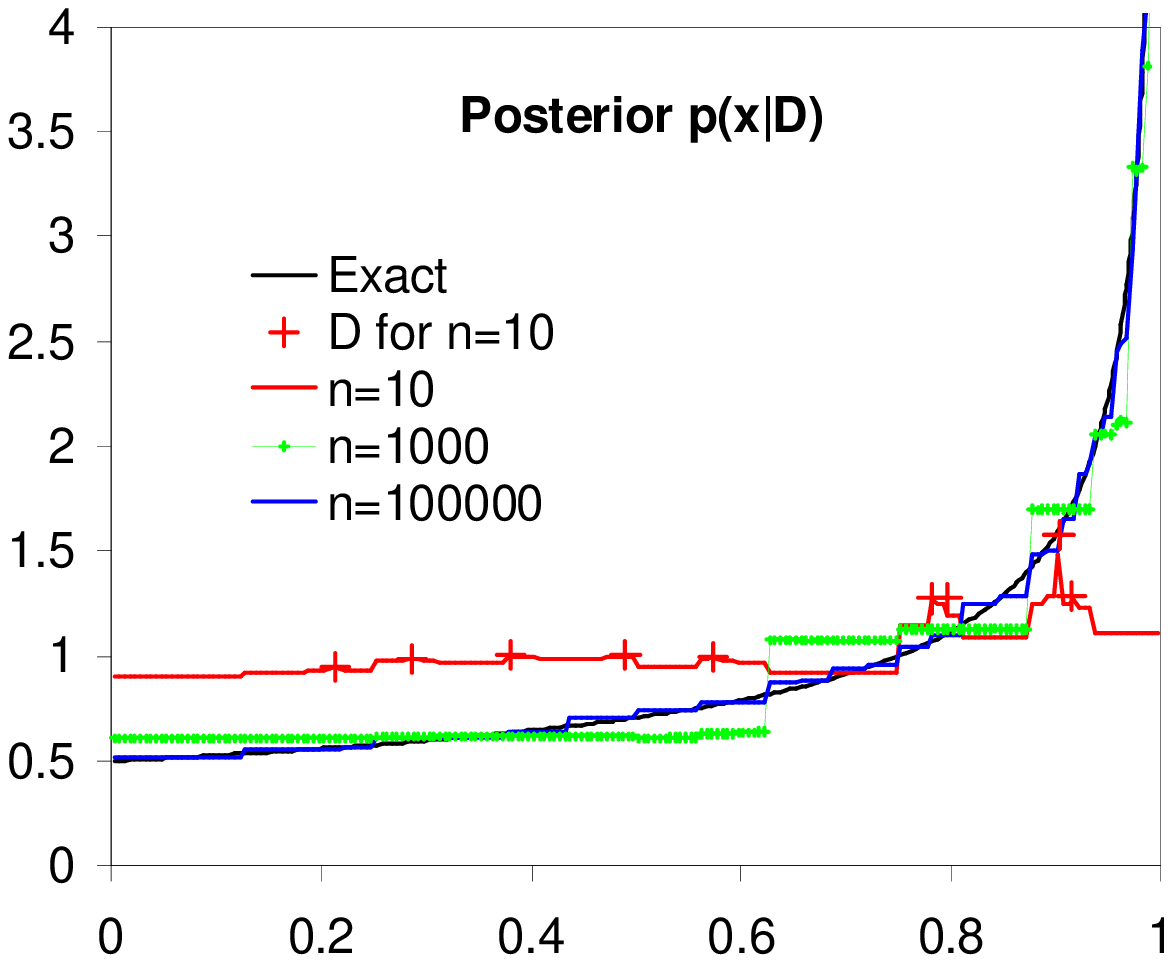}{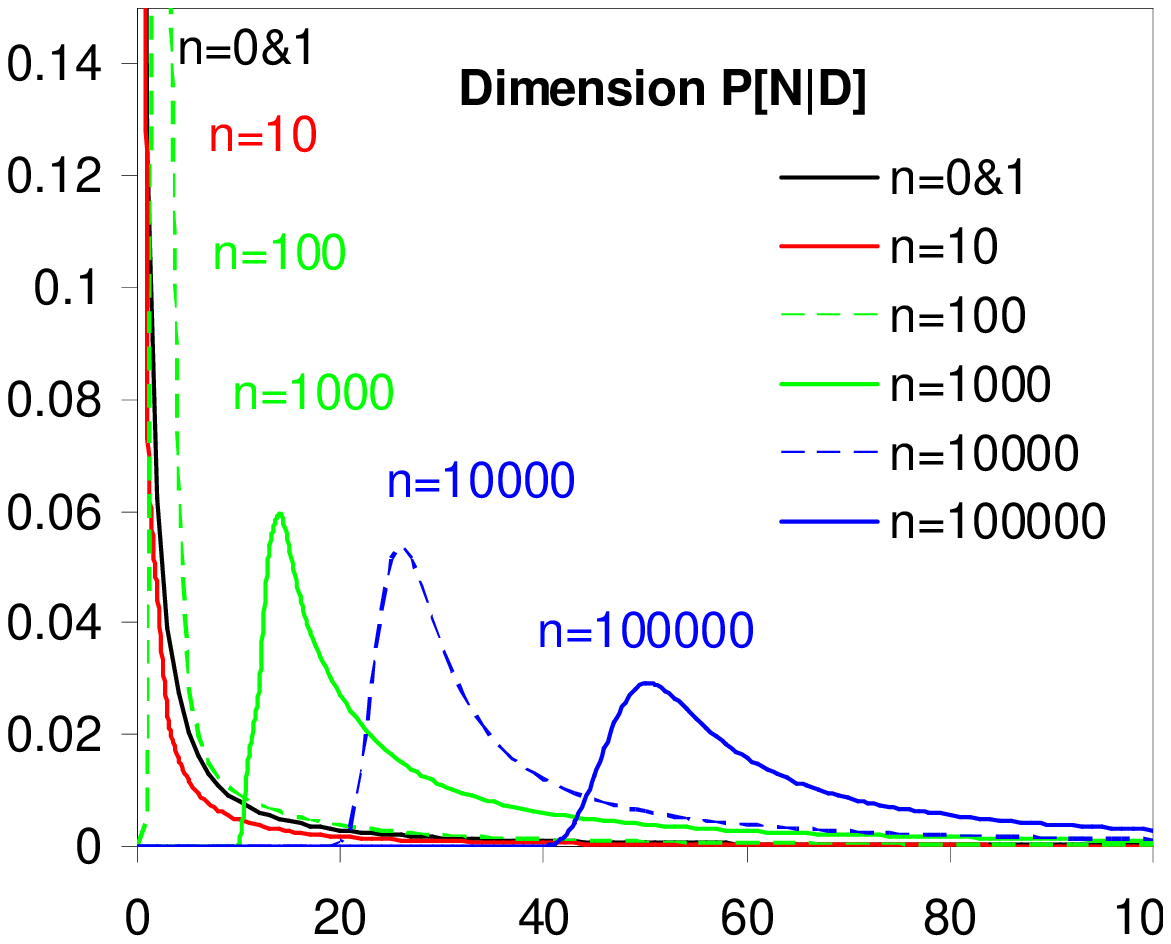}%
{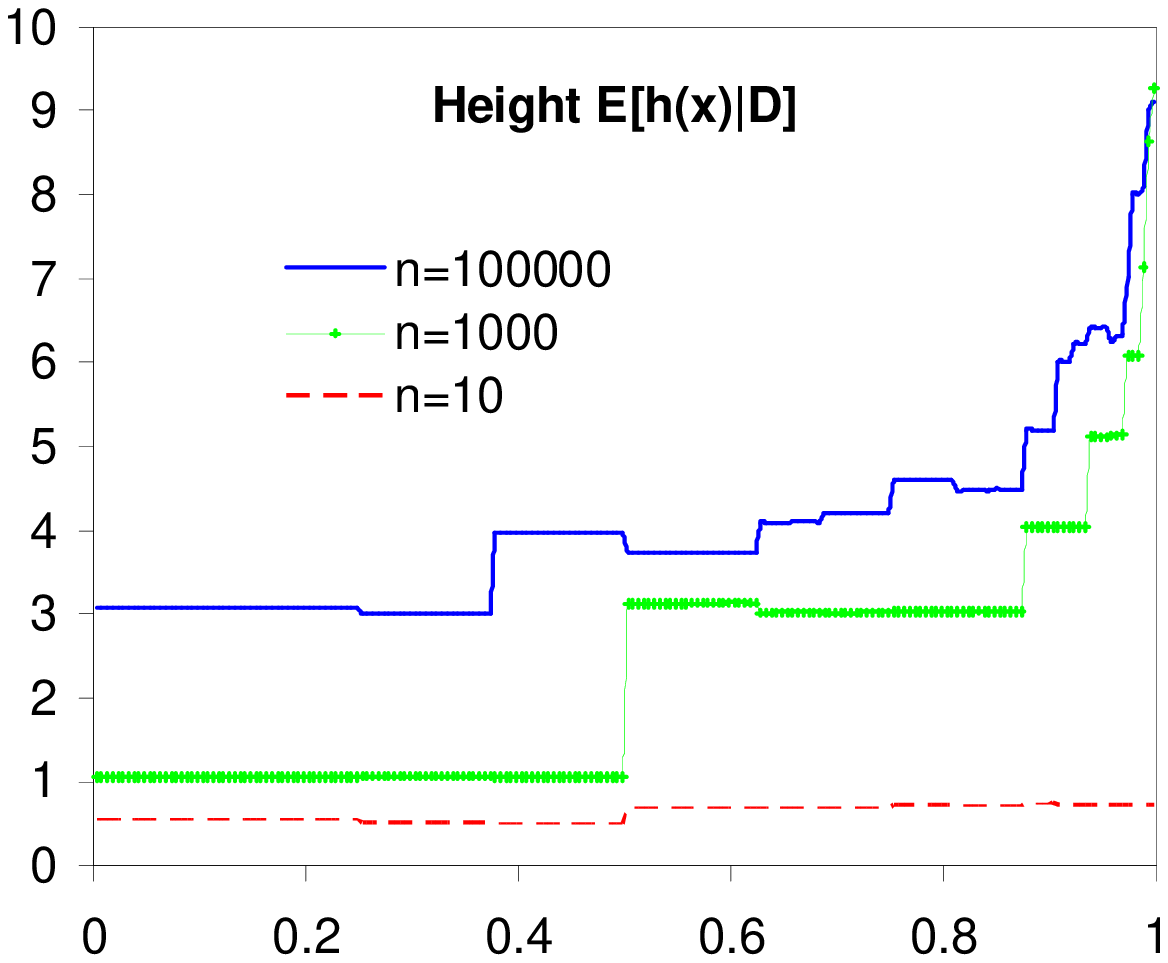}{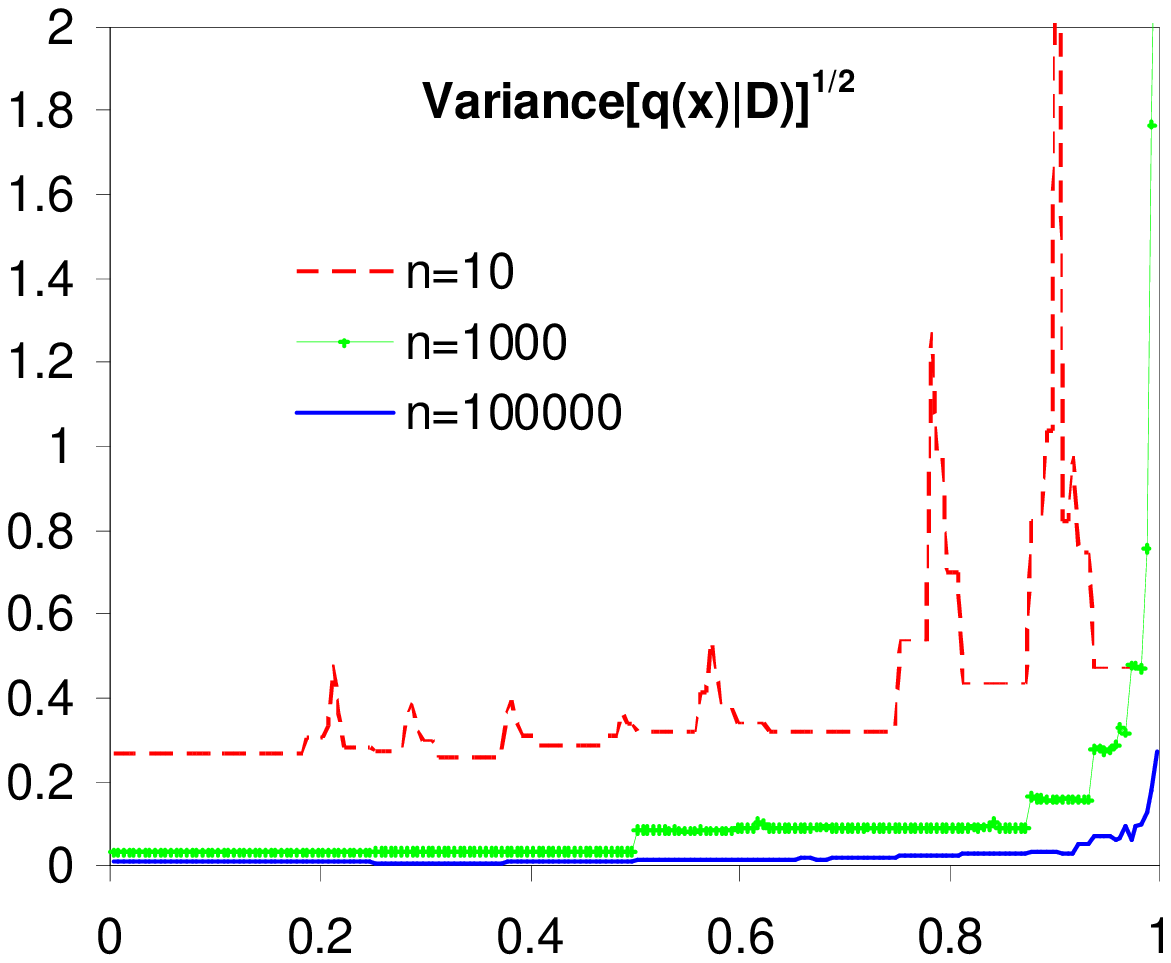}{figSingular}%
{BayesTree() results for the {\bf Singular}
distribution $\dot q(x)=2/\sqrt{1-x}$, prototype for a proper
singular distribution.}

\fourgraphs{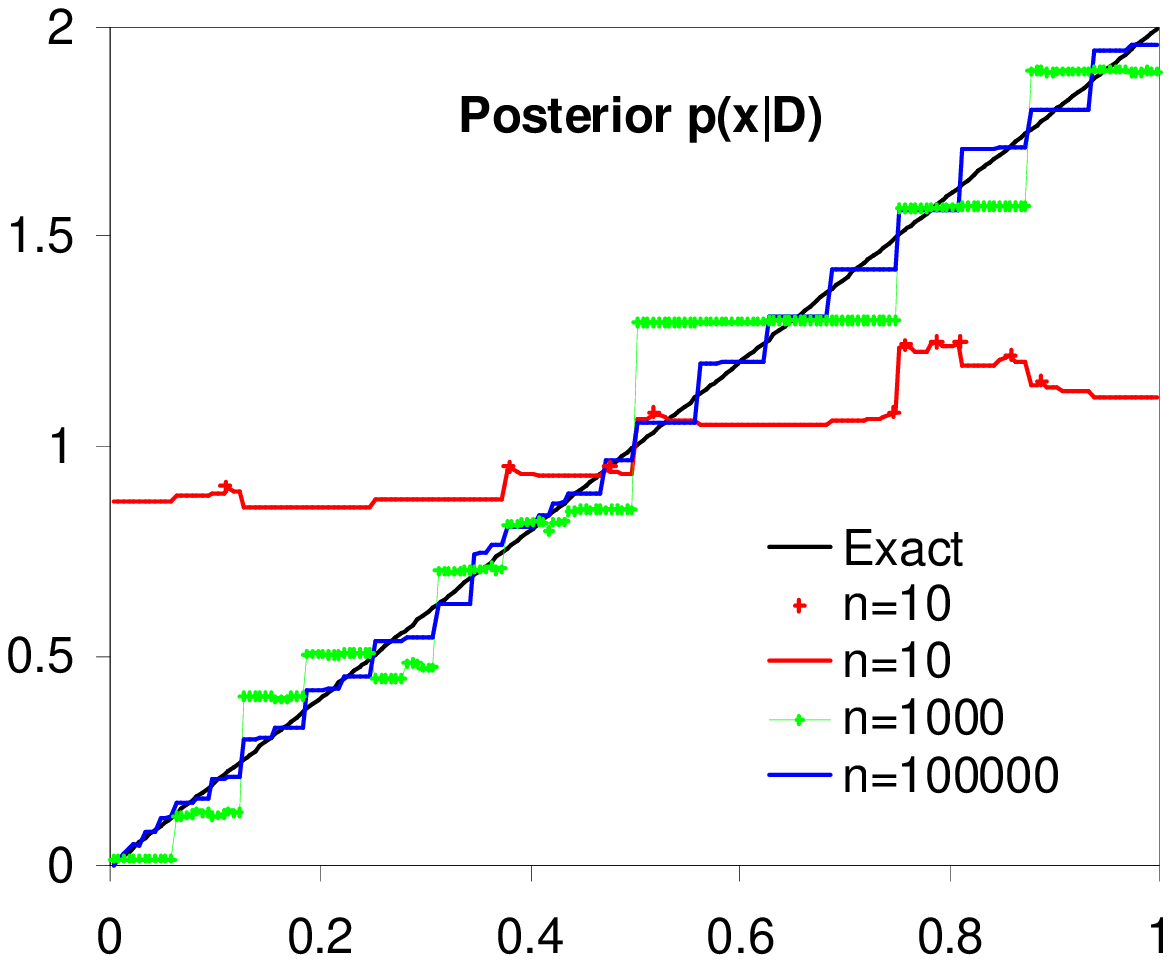}{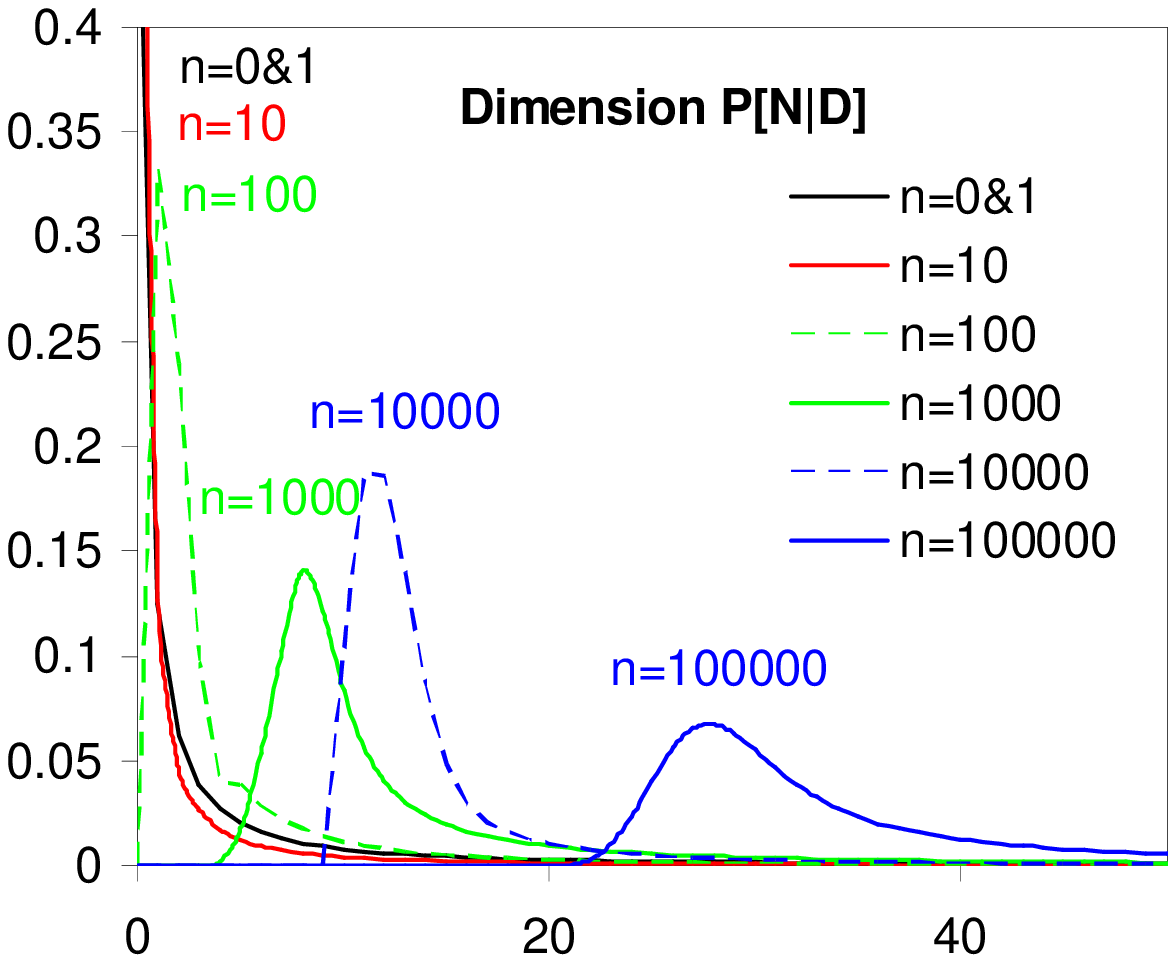}%
{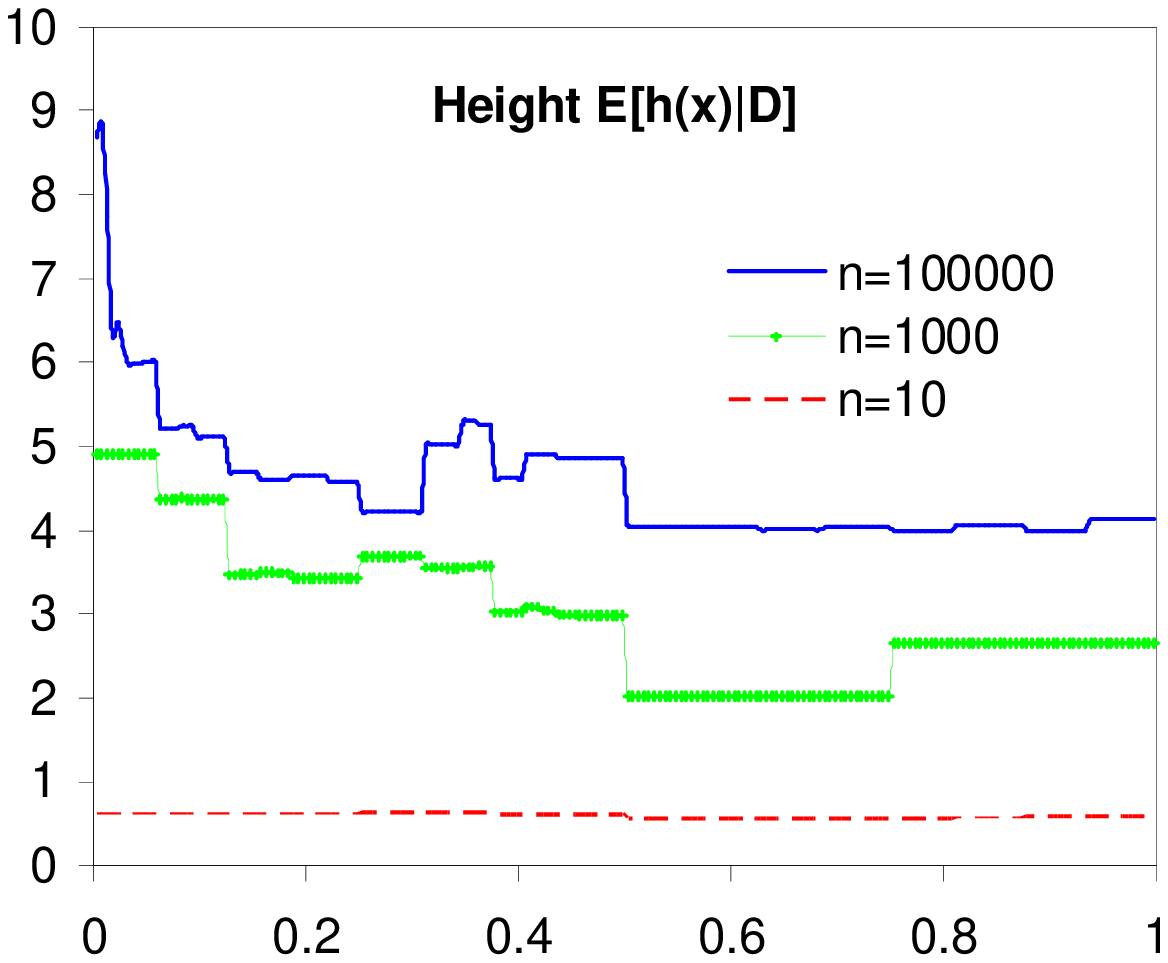}{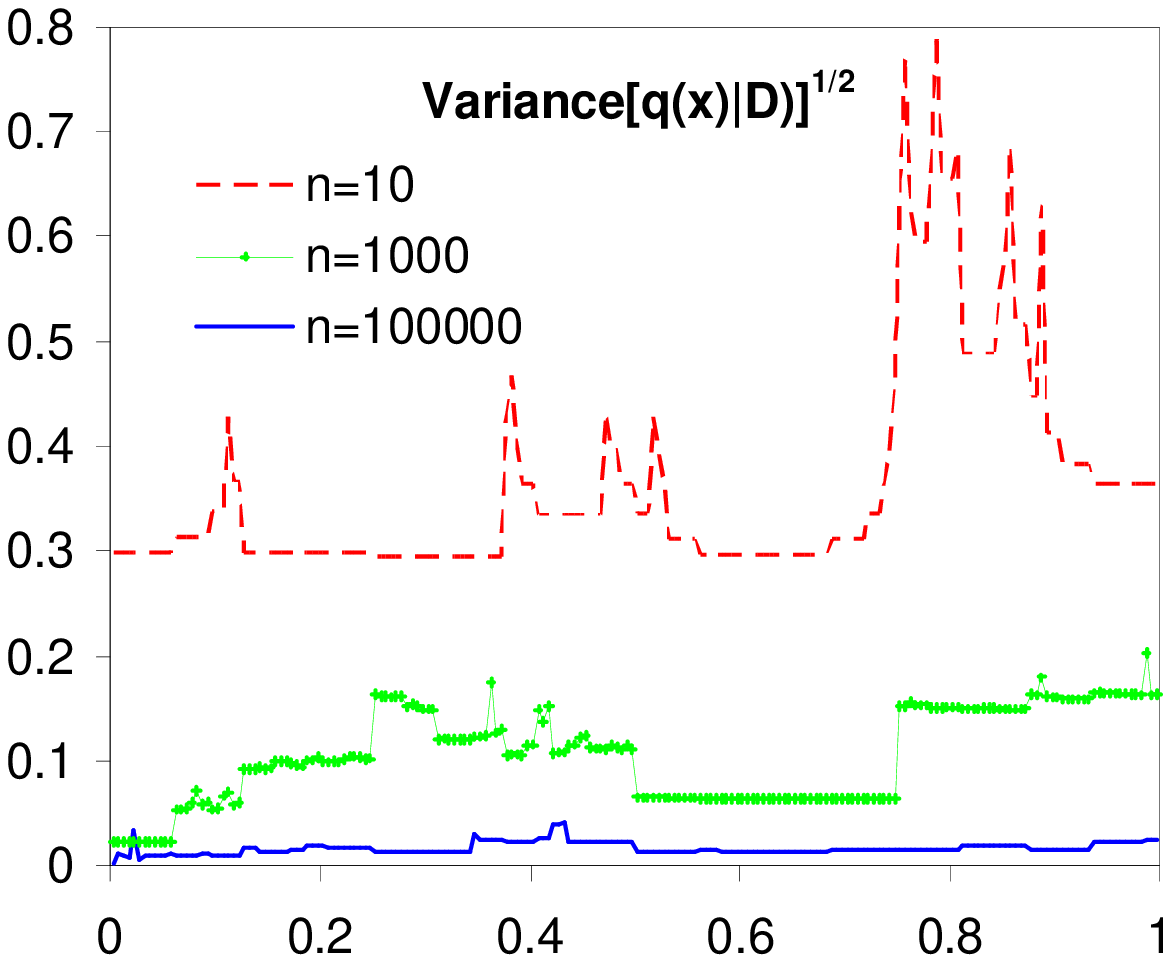}{figLinear}%
{BayesTree() results for the {\bf Linear}
distribution $\dot q(x)=2x$, prototype for a continuous function.}

\fourgraphs{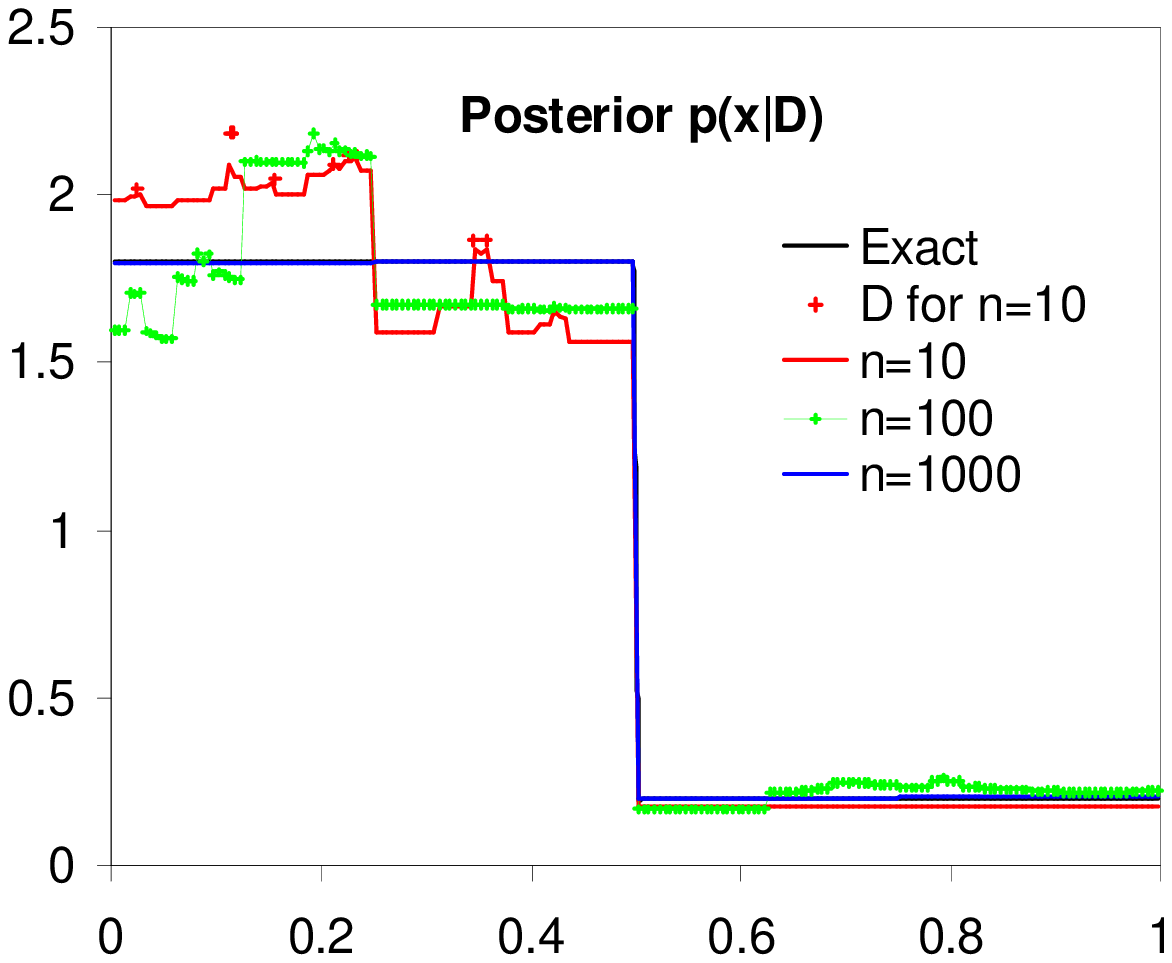}{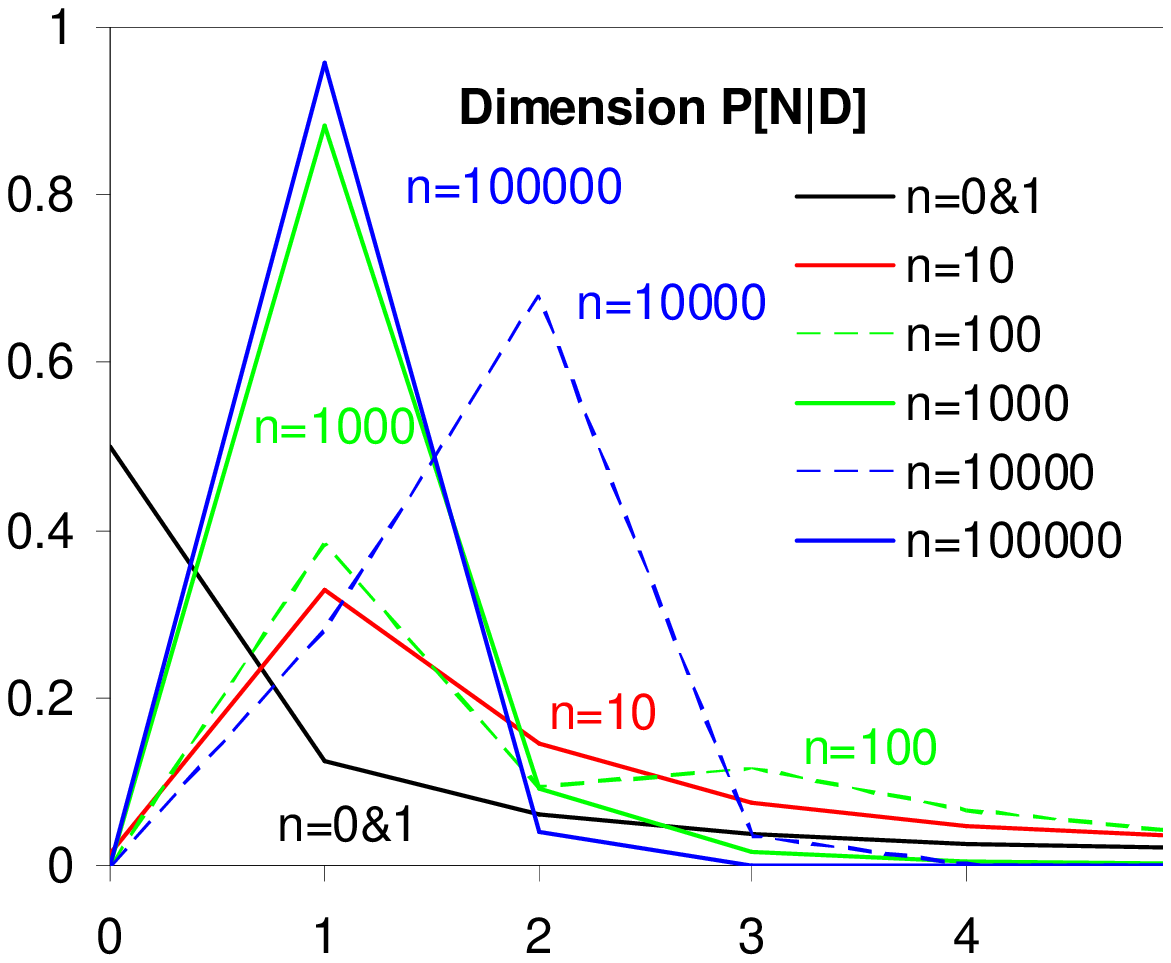}%
{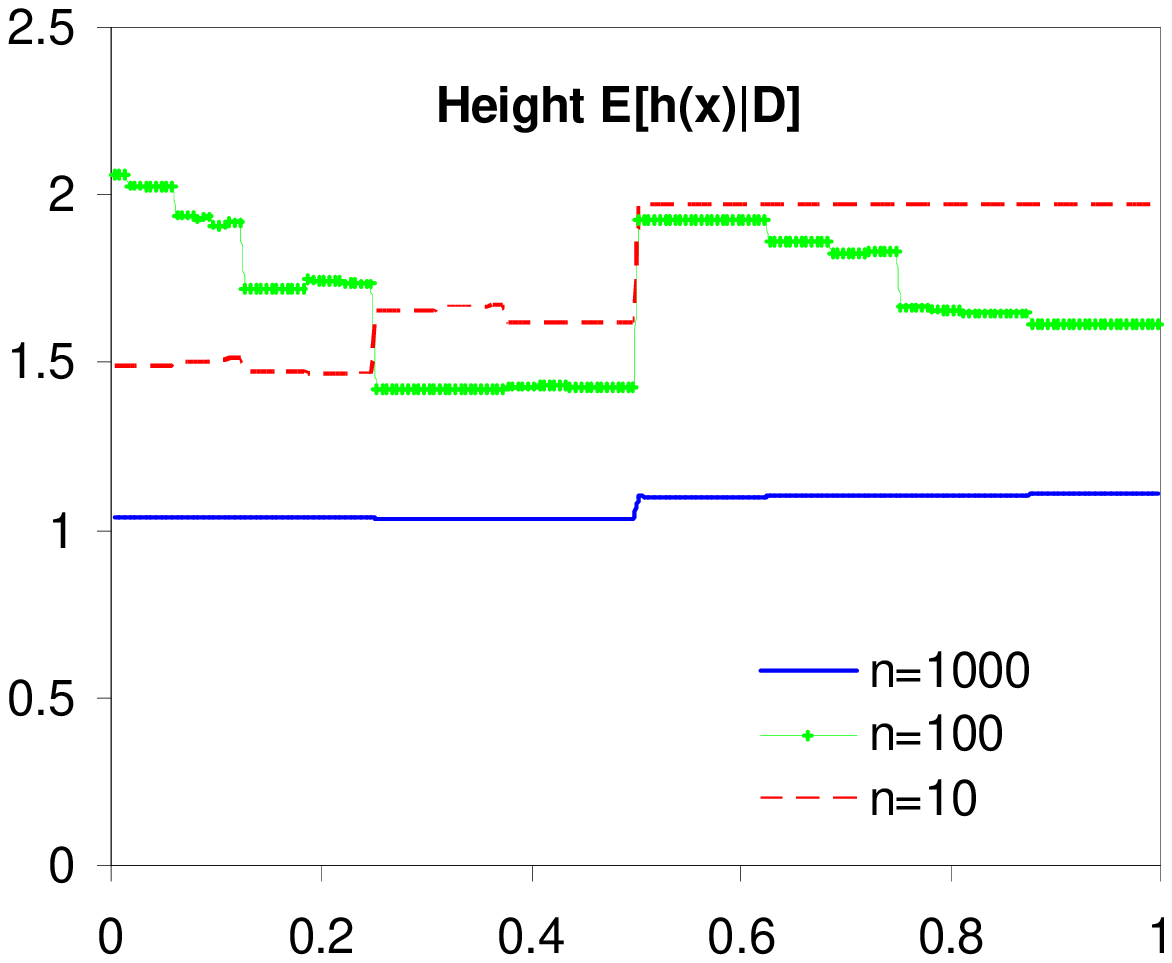}{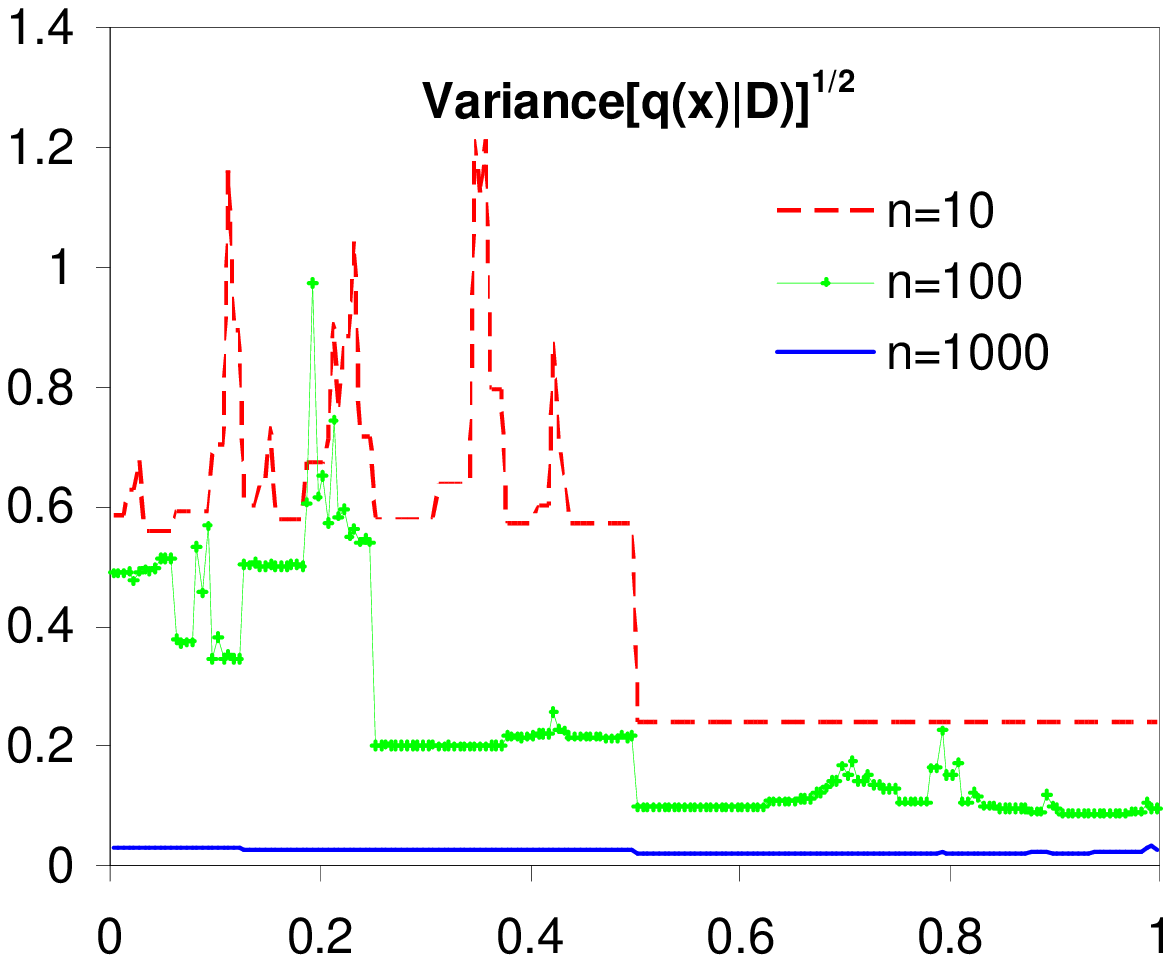}{figJump05}%
{BayesTree() results for the {\bf Jump-at-1/2} distribution $\dot
q(x)=9/5$ for $x<1/2$ and $q(x)=1/5$ for $x\geq 1/2$, prototype
for a piecewise constant function with a finite Bayes tree.}

\fourgraphs{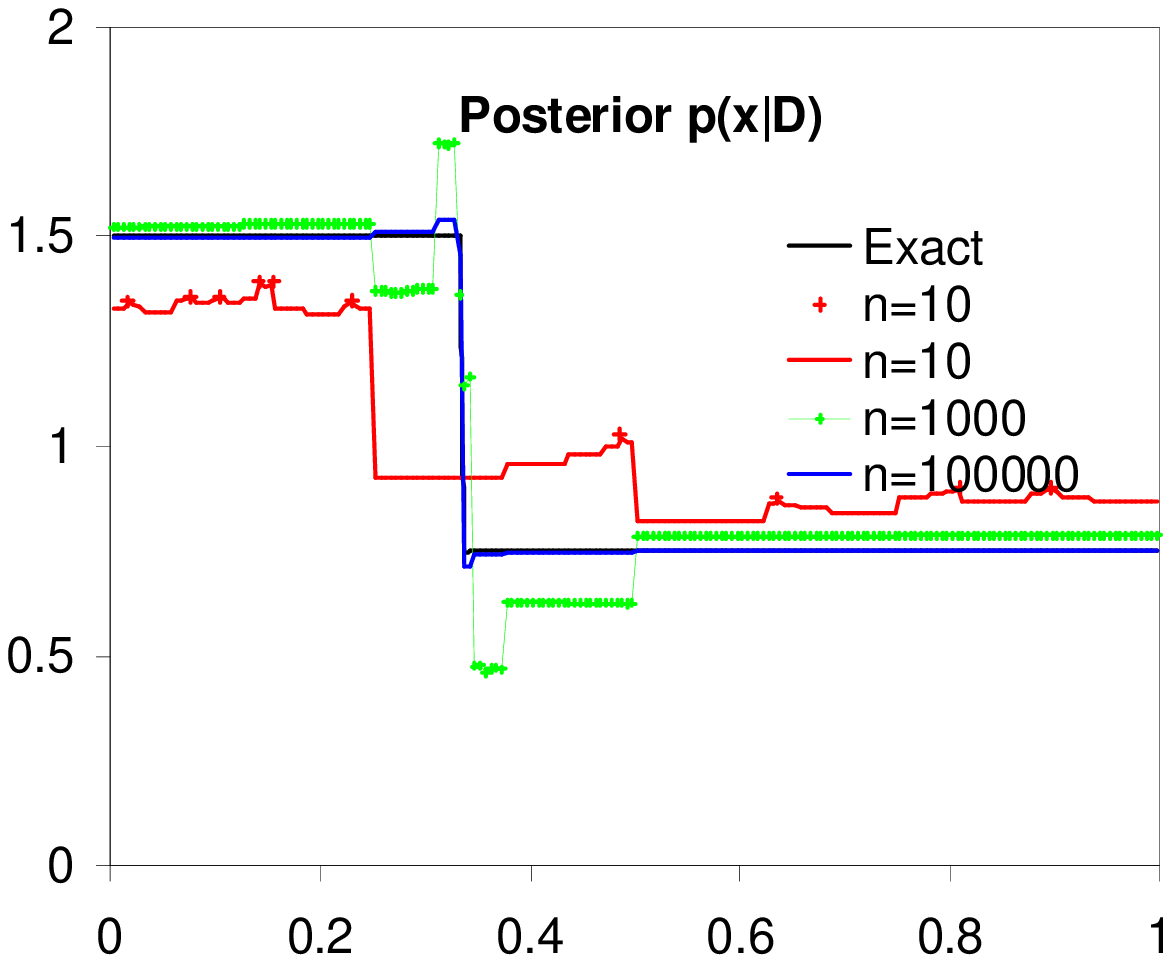}{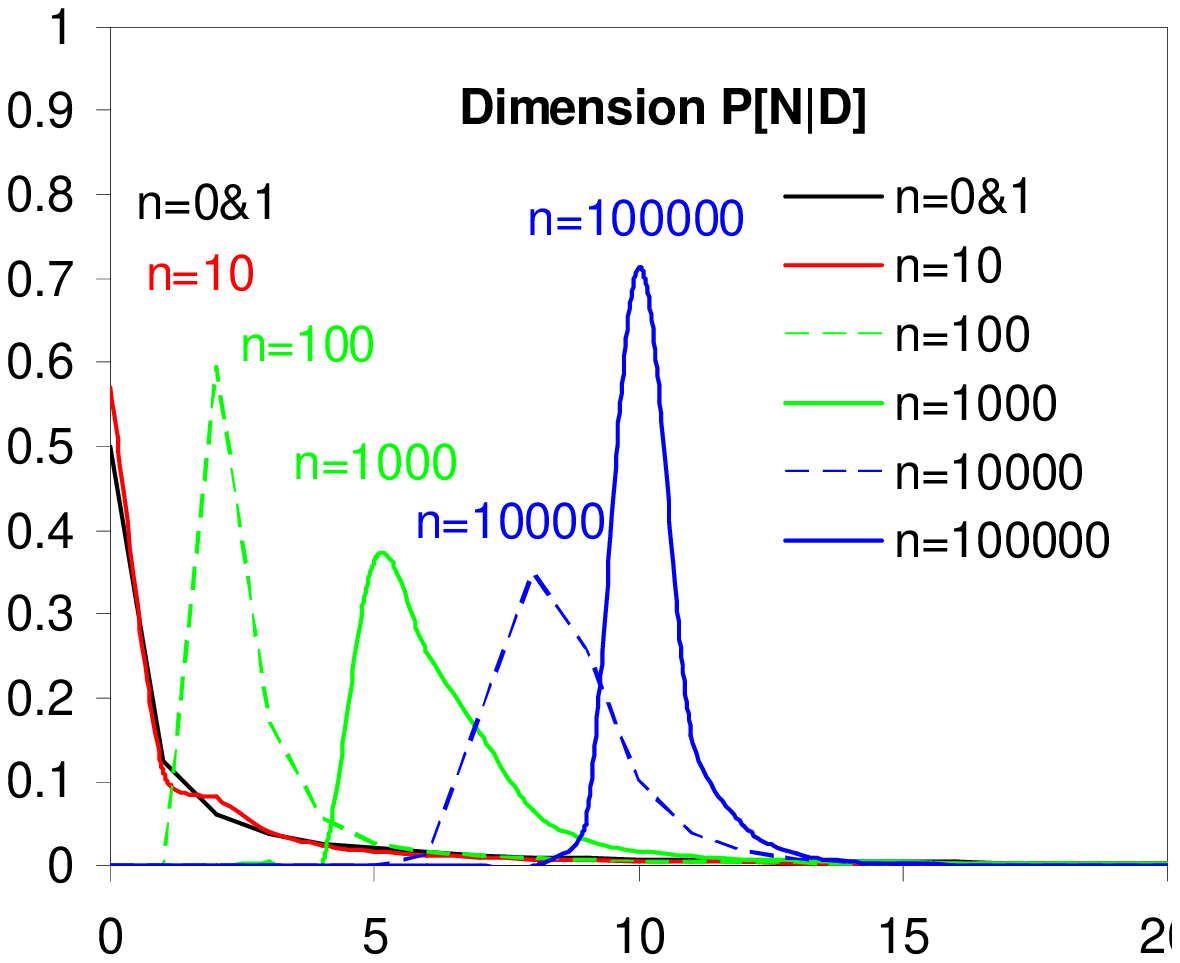}%
{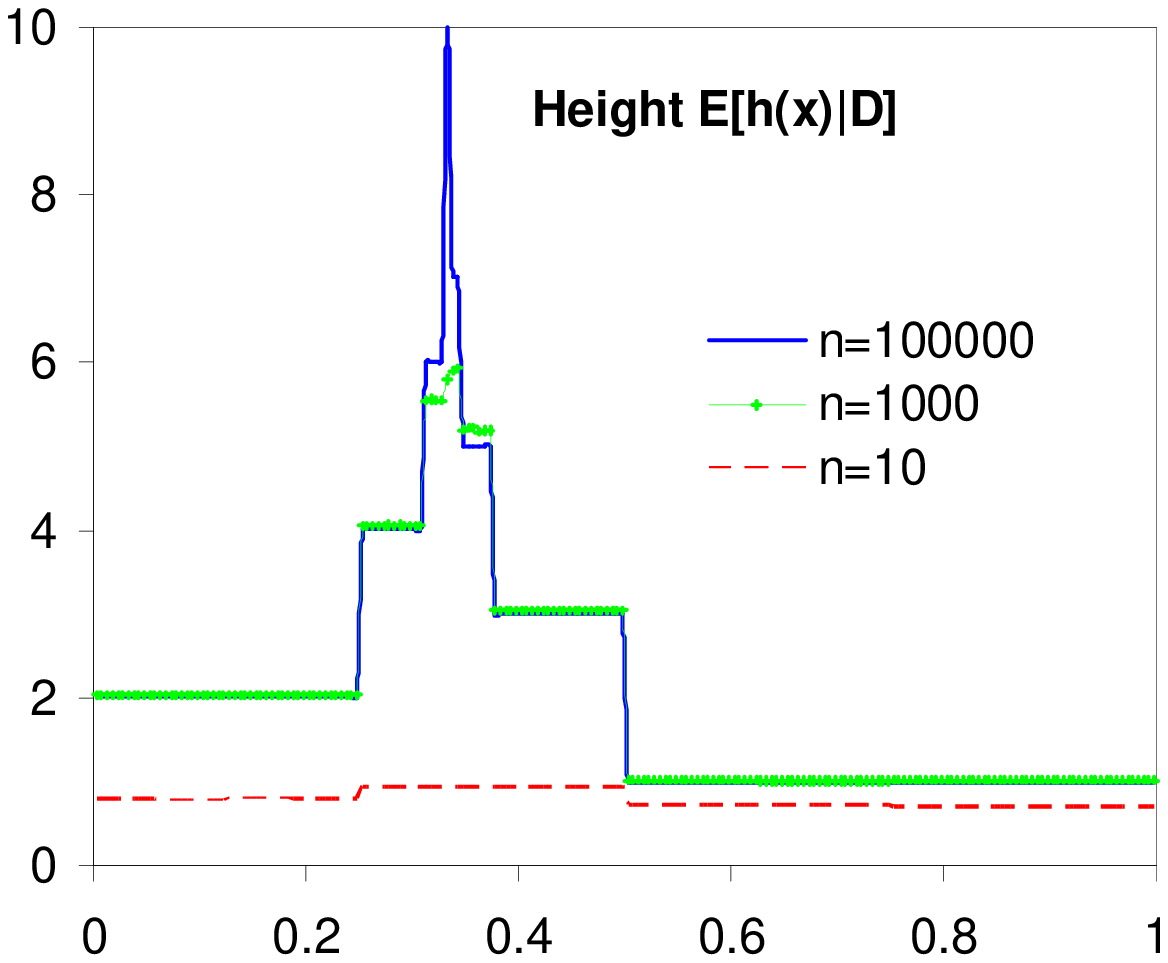}{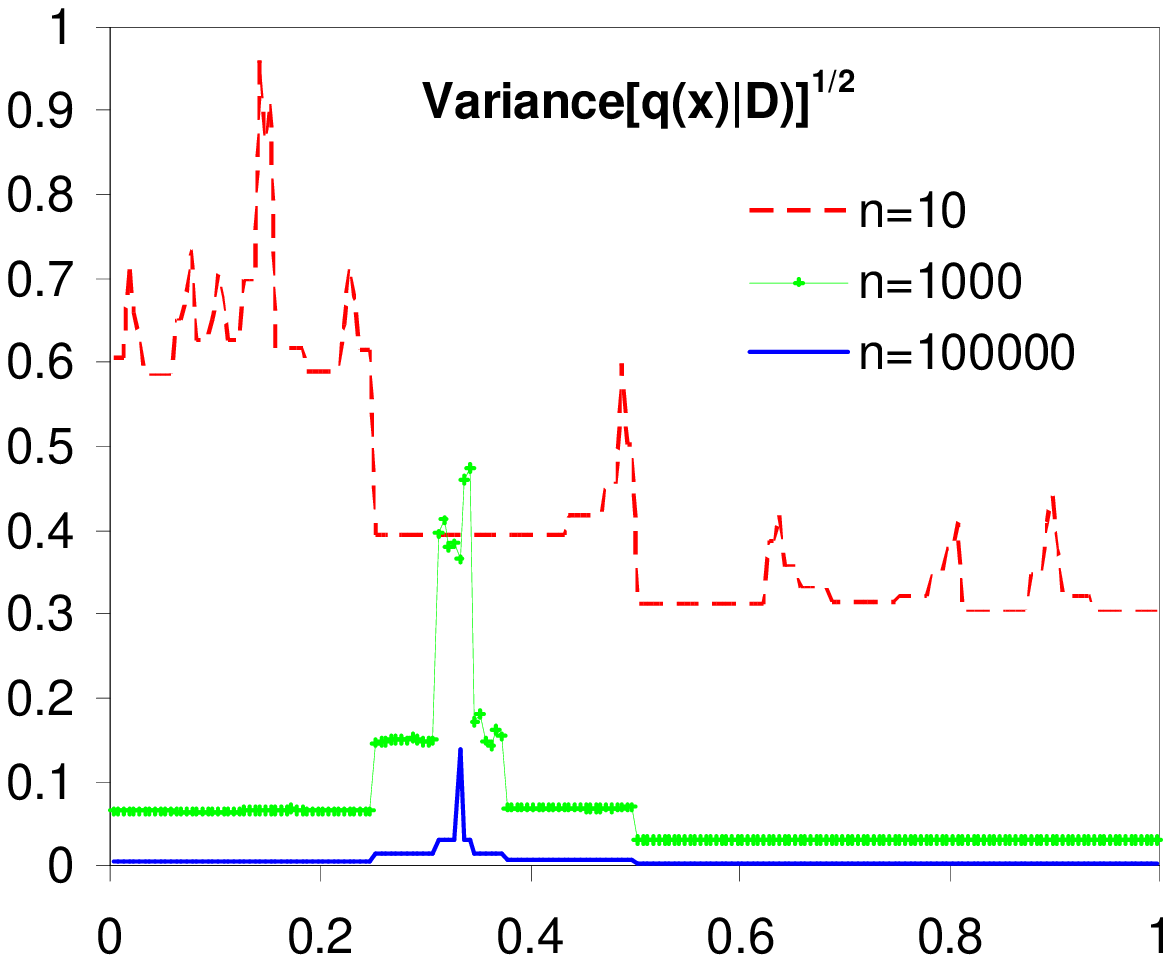}{figJump}%
{BayesTree() results for the {\bf Jump-at-1/3} distribution $\dot
q(x)=2/3$ for $x<{\bf 1/3}$ and $q(x)=1/3$ for $x\geq 1/3$,
prototype for a piecewise constant function with an infinite Bayes
tree.}

\paradot{Graphs and examples}
To get further insight into the behavior of our model, we
numerically investigated some example distributions $\dot q()$. We
have chosen elementary functions, which can be regarded as
prototypes for more realistic functions. They include the Beta,
linear, a singular, and piecewise constant distributions with finite
and infinite Bayes trees. These examples on $[0,1)$ also shed light
on the other spaces discussed in Section \ref{secTMM}, since they
are isomorphic. The posteriors, model dimensions, tree heights, and
variances are plotted in Figures \ref{figBeta}--\ref{figJump} for
random samples $D$ of sizes $n=10^0,...,10^5$. We first discuss
observations common to all sampling distributions, thereafter
specific aspects. All experiments were performed with split
probability $s=\fr12$ and uniform distribution $\a=1$.

\paradot{General observations}
The posteriors $p(x|D)$ clearly converge for $n\to\infty$ to the
true distribution $\dot q()$, accompanied by a (necessary)
moderate growth of the effective dimension (except for
Jump-at-1/2). For $n=10$ we show the data points. It is visible
how each data point pulls the posterior up, as it should be (``one
sample seldom comes alone'').

Compare this to an empirical bin model with $N$ bins. Since each
bin contains $O(n/N)$ data points, the frequency estimate
$n_{bin}/n$ of $\dot q(bin)$ has accuracy $O(\sqrt{N/n})$. The
minimal error when approximating a continuous function by a
piecewise constant function with bin size $1/N$ is $O(1/N)$, so the
estimate has total error $\max\{O(\sqrt{N/n}),O(1/N)\}$ with
minimum $O(n^{-1/3})$ at\footnote{Sometimes heuristic $N=\sqrt{n}$ is
proposed, which makes no sense.} $N=n^{1/3}$. This is nicely
consistent with our model dimension. Look at the maxima of dimension
distribution $P[N|D]$ or count the number of significant jumps in the
posterior $p(x|D)$.

The figures also show that the posterior variances Var[$q(x)|D$]
converge to zero for $n\to\infty$, but diverge when $x$ tends to
a point in $D$, consistent with the theoretical analysis of
multi-points in Section \ref{secIT}.
The expected tree height $E[h(x)|D]$ at $x$ correctly reflects the
local needs for (non)splits.

\paradot{Specific observations}
{\it Beta:} The Beta distribution Beta$(x|\alpha,\beta)\propto
x^{\alpha-1}(1-x)^{\beta-1}$ is prototypical for a smooth unimodal
distribution. Apart from local jitter, the simplest model
consistent with data size $n=10$ is essentially a Jump-at-1/2
function (see below). The tree height slowly increases with $n$
with a dip around the ``flat'' maximum of the Beta, since a
constant approximation works well there.

{\it Singular:} We used the distribution $\dot q(x)=2/\sqrt{1-x}$ as a
prototype for a proper singular distribution. The tree height
is necessarily larger near the singularity at $x=1$.

{\it Linear:} Once continuously differentiable functions are locally
linear, so the linear distribution $\dot q(x)=2x$ serves as a
prototype for them. The better approximation of $p(x|D)$ near 0
versus near 1, accompanied by a higher tree, is remarkable. First,
there are fewer data points near 0 to warrant this, and second, the
region is less interesting, since more samples are at 1. So we
expected quite the opposite behavior. We currently have no
explanation for this phenomenon.

{\it Jump-at-1/2:} Also illustrative are distributions with finite
Bayes tree, i.e.\ piecewise constant functions with
discontinuities only at binary fractions. We consider the
prototype that jumps at $x=1/2$. All quantities converge rapidly.
We see that model dimension and tree height stay finite in this
case, as they should. Both converge to the minimal consistent value 1.
The variance in the left and right half of
$[0,1)$ is roughly proportional to $\dot q$ therein.

{\it Jump-at-1/3:} A jump at a non-binary fraction cannot be
modeled with a finite tree. Convergence is slower than for
Jump-at-1/2, but faster than for the other examples,
which makes sense since only one branch of the tree has to grow to
infinity. This shows up in a slower increase of dimension, a
converging height function with singularity at 1/3, and a
narrowing spike in the variance.

\section{Discussion}\label{secDisc}

We presented a Bayesian model on infinite trees, where we split a
node into two subtrees with some probability, and assigned a
Beta distributed probability to each subtree. %
We were primarily interested in the case of zero prior knowledge.
In this case, scale invariance and symmetry should be preserved. %
Scale invariance requires the parameters to be the same in each
node and symmetry requires a symmetric Beta, leaving one
splitting probability $s$ and one Beta parameter $\beta$ as
adjustable parameters. %
We devised closed form expressions for various inferential
quantities of interest at the data separation level, which led to
an exact algorithm with runtime essentially linear in the data
size.

We introduced and studied this two-parameter tree-model class. The
most interesting case of splitting probability $s=\fr12$ and uniform
prior over subtrees $\beta=1$ has been studied in more detail. The
theoretical and numerical model behavior was very reasonable, e.g.\
consistency (no underfitting) and low finite effective dimension (no
overfitting). Higher moments can be made finite by smaller $s$ or
larger $\beta$.

There are various natural generalizations of our model. %
The splitting probability $s$ and Beta parameter $\beta$ could be
made dependent on the node of the tree, which allows incorporating
prior knowledge. $k$-ary trees could be allowed with Beta
generalized to Dirichlet distributions. %
Non-symmetric partitions are straightforward to implement by replacing
all $\delta(q_z-\fr12)$ with $\delta(q_z-|\G_z|/|\G_{z_{1:l-1}}|)$, and
possibly using non-symmetric Betas. %
The expected entropy can also be computed by allowing fractional
counts $n_z$ and noting that $x\ln x = {d\over dx}
x^\alpha|_{\alpha=1}$ \cite{Hutter:01xentropy,Hutter:05mifs}. %
A sort of maximum a posteriori (MAP) tree skeleton can also easily be
read off from (\ref{eqEvDens}). A node $\G_z$ in the MAP-like tree
is a leaf iff $s{p_{z0}(D_{z0})p_{z1}(D_{z1})\over
w(n_{z0},n_{z1})}<u$. %
A challenge is to generalize the model from piecewise constant to
piecewise linear continuous functions, at least for
$\G=[0,1)$. Independence of subtrees no longer holds,
which was key in our analysis.

If $\G$ is not already a tree or binary string, but an interval, a
major problem of Polya trees and our tree model are partition
artifacts in the estimated density. Numerically but unlikely
analytically it is possible to average over boundary locations like
in \cite{Paddock:03} and smooth out discontinuities. Interestingly
for flat bin estimation, analytical averaging {\em is} possible via
dynamic programming \cite{Endres:05}.


\begin{small}

\end{small}


\begin{thebibliography}{ABCD}\parskip=0ex

\bibitem[Bis06]{Bishop:06}
C.~M. Bishop.
\newblock {\em Pattern Recognition and Machine Learning}.
\newblock Springer, 2006.

\bibitem[BM98]{Borovkov:98}
A.~A. Borovkov and A.~Moullagaliev.
\newblock {\em Mathematical Statistics}.
\newblock Gordon \& Breach, 1998.

\bibitem[DLR77]{Dempster:77}
A.~P. Dempster, N.~Laird, and D.~Rubin.
\newblock Maximum likelihood estimation for incomplete data via the {EM}
  algorithm.
\newblock {\em Journal of the Royal Statistical Society}, Series B 39:1--38,
  1977.

\bibitem[EF05]{Endres:05}
D.~Endres and P.~F{\"o}ldi{\'a}k.
\newblock Bayesian bin distribution inference and mutual information.
\newblock {\em IEEE Transactions on Information Theory}, 51(11):3766--3779,
  2005.

\bibitem[EW95]{Escobar:95}
M.~Escobar and M.~West.
\newblock Bayesian density estimation and inference using mixtures.
\newblock {\em Journal of the American Statistical Association}, 90:577--588,
  1995.

\bibitem[Fab64]{Fabius:64}
J.~Fabius.
\newblock Asymptotic behavior of {B}ayes estimates.
\newblock {\em Annals of Mathematical Statistics}, 35:846--856, 1964.

\bibitem[Fer73]{Ferguson:73}
T.~S. Ferguson.
\newblock On the mathematical foundations of theoretical statistics.
\newblock {\em Annals of Statistics}, 1(2):209--230, 1973.

\bibitem[GM03]{Gray:03}
A.~G. Gray and A.~W. Moore.
\newblock Nonparametric density estimation: Toward computational tractability.
\newblock In {\em SIAM International Conf. on Data Mining}, volume~3, 2003.

\bibitem[Goo83]{Good:83}
I.~J. Good.
\newblock Explicativity, corroboration, and the relative odds of hypotheses.
\newblock In {\em Good thinking: The Foundations of Probability and its
  applications}. University of Minnesota Press, Minneapolis, MN, 1983.

\bibitem[Goo84]{Good:84}
I.~J. Good.
\newblock The best explicatum for weight of evidence.
\newblock {\em Journal of Statistical Computation and Simulation}, 19:294--299,
  1984.

\bibitem[Hut02]{Hutter:01xentropy}
M.~Hutter.
\newblock Distribution of mutual information.
\newblock In {\em Advances in Neural Information Processing Systems 14
  ({NIPS'01})}, pages 399--406, Cambridge, MA, 2002. MIT Press.

\bibitem[Hut05a]{Hutter:05bayestree}
M.~Hutter.
\newblock Fast non-parametric {B}ayesian inference on infinite trees.
\newblock In {\em Proc. 10th International Conf. on Artificial Intelligence and
  Statistics ({AISTATS'05})}, pages 144--151. Society for Artificial
  Intelligence and Statistics, 2005.

\bibitem[Hut05b]{Hutter:04uaibook}
M.~Hutter.
\newblock {\em Universal Artificial Intelligence: Sequential Decisions based on
  Algorithmic Probability}.
\newblock Springer, Berlin, 2005.
\newblock 300 pages, http://www.hutter1.net/ai/uaibook.htm.

\bibitem[Hut07]{Hutter:07btcode}
M.~Hutter.
\newblock Additional material to article, 2007.
\newblock \\ http://www.hutter1.net/official/bib.htm\#bayestreex.

\bibitem[HZ05]{Hutter:05mifs}
M.~Hutter and M.~Zaffalon.
\newblock Distribution of mutual information from complete and incomplete data.
\newblock {\em Computational Statistics \& Data Analysis}, 48(3):633--657,
  2005.

\bibitem[Jay03]{Jaynes:03}
E.~T. Jaynes.
\newblock {\em Probability Theory: The Logic of Science}.
\newblock Cambridge University Press, Cambridge, MA, 2003.

\bibitem[Jef35]{Jeffreys:35}
H.~Jeffreys.
\newblock Some tests of significance, treated by the theory of probability.
\newblock In {\em Proc. Cambridge Philosophical Society}, volume~31, pages
  203--222, 1935.

\bibitem[KF98]{Koller:98}
D.~Koller and R.~Fratkina.
\newblock Using learning for approximation in stochastic processes.
\newblock In {\em Proc. 15th International Conference on Machine Learning
  ({ICML'98})}, pages 287--295, 1998.

\bibitem[KK97]{Kozlov:97}
A.~V. Kozlov and D.~Koller.
\newblock Nonuniform dynamic discretization in hybrid networks.
\newblock In {\em Proc. 13th Conf. on Uncertainty in Artificial Intelligence
  ({UAI'97})}, pages 314--325, 1997.

\bibitem[KM07]{Kontkanen:07}
P.~Kontkanen and P.~Myllymäki.
\newblock {MDL} histogram density estimation.
\newblock In {\em Proc. 11th International Conf. on Artificial Intelligence and
  Statistics ({AISTATS'07})}. Society for Artificial Intelligence and
  Statistics, 2007.

\bibitem[Lav92]{Lavine:92}
M.~Lavine.
\newblock Some aspects of {P}olya tree distributions for statistical modelling.
\newblock {\em Annals of Statistics}, 20:1222--1235, 1992.

\bibitem[Lav94]{Lavine:94}
M.~Lavine.
\newblock More aspects of {P}olya tree distributions for statistical modelling.
\newblock {\em Annals of Statistics}, 22:1161--1176, 1994.

\bibitem[Lem03]{Lemm:03}
J.~C. Lemm.
\newblock {\em Bayesian Field Theory and Approximate Symmetries}.
\newblock Johns Hopkins University Press, 2003.

\bibitem[LLW07]{Liu:07}
H.~Liu, J.~Lafferty, and L.~Wasserman.
\newblock Sparse nonparametric density estimation in high dimensions using the
  rodeo.
\newblock In {\em Proc. 11th International Conf. on Artificial Intelligence and
  Statistics ({AISTATS'07})}. Society for Artificial Intelligence and
  Statistics, 2007.

\bibitem[Mac03]{MacKay:03}
D.~J.~C. MacKay.
\newblock {\em Information theory, inference and learning algorithms}.
\newblock Cambridge University Press, Cambridge, MA, 2003.

\bibitem[PRLW03]{Paddock:03}
S.~M. Paddock, F.~Ruggeri, M.~Lavine, and M.~West.
\newblock Randomised {P}olya tree models for nonparametric {B}ayesian
  inference.
\newblock {\em Statistica Sinica}, 13(2):443--460, 2003.

\bibitem[PW02]{Petrone:02}
S.~Petrone and L.~Wasserman.
\newblock Consistency of {B}ernstein polynomial posteriors.
\newblock {\em Journal of the Royal Statistical Society}, B 64:79--100, 2002.

\end{thebibliography}
\end{document}